\documentclass[11pt]{amsart}

\usepackage{amsmath,amsfonts,amssymb,amsthm}
\usepackage[cmtip,all]{xy}
\usepackage[numbers,square,sort&compress]{natbib}
\usepackage{subfig}
\usepackage{xspace}
\usepackage{a4wide}
\usepackage{tikz}
\usepackage{tkz-graph}
\usetikzlibrary{decorations.markings}
\usetikzlibrary{arrows}

\graphicspath{{./}{./figs/}{../../figs/}}

\def\N{\mathbb{N}}
\def\Z{\mathbb{Z}}

\def\simkappa{\!\!\sim_{\kappa}}
\def\tsimkappa{\sim_{\kappa}}

\def\simalpha{\!\!\sim_{\alpha}}
\def\tsimalpha{\sim_{\alpha}}
\def\SYmod{S_\Gamma\!/\simalpha}
\def\sds{{SDS}\xspace}
\def\sdss{{SDS}s\xspace}
\def\Acyc{\mathrm{Acyc}}
\def\acyc{\alpha}
\DeclareMathOperator{\C}{C}

\def\<{\langle}
\def\>{\rangle}
\def\card#1{|#1|}
\DeclareMathOperator{\Line}{Line}
\DeclareMathOperator{\Circle}{Circle}
\def\vset{\mathrm{v}}
\def\eset{\mathrm{e}}
\def\eg{e.g.\xspace}
\def\ie{i.e.\xspace}
\def\nkeq{\nsim_\kappa}
\def\varepsilon{\eta_e}
\def\ieext{\mathcal{I}^*_e}

\newtheorem{defn}{Definition}[section]
\newtheorem{prop}[defn]{Proposition}

\newtheorem{thm}[defn]{Theorem}
\newtheorem{cor}[defn]{Corollary}

\newtheorem{ex}[defn]{Example}

\begin{document}

\title{Posets from Admissible Coxeter Sequences}

\author{Matthew Macauley} 
\author{Henning S.~Mortveit} 

\begin{abstract}
  We study the equivalence relation on the set of acyclic orientations
  of an undirected graph $\Gamma$ generated by source-to-sink
  conversions. These conversions arise in the contexts of admissible
  sequences in Coxeter theory, quiver representations, and
  asynchronous graph dynamical systems. To each equivalence class we
  associate a poset, characterize combinatorial properties of these
  posets, and in turn, the admissible sequences. This allows us to
  construct an explicit bijection from the equivalence classes over
  $\Gamma$ to those over $\Gamma'$ and $\Gamma''$, the graphs obtained
  from $\Gamma$ by edge deletion and edge contraction of a fixed
  cycle-edge, respectively.  This bijection yields quick and elegant
  proofs of two non-trivial results: $(i)$ A complete combinatorial
  invariant of the equivalence classes, and $(ii)$ a solution to the
  conjugacy problem of Coxeter elements for simply-laced Coxeter
  groups. The latter was recently proven by H.~Eriksson and
  K.~Eriksson using a much different approach.
\end{abstract}



\maketitle


\section{Overview.}
\label{sec:overview}

Let $O_\Gamma$ be an acyclic orientation of the undirected graph
$\Gamma$. A cyclic $1$-shift (left) of a linear extension $\pi$ of
$O_\Gamma$ corresponds to converting a source (the element $\pi_1$) of
$O_\Gamma$ into a sink, and this gives rise to an equivalence relation
on $\Acyc(\Gamma)$ denoted by $\tsimkappa$. We let $\kappa(\Gamma)$
denote the number of equivalence classes in $\Acyc(\Gamma)$ under
$\tsimkappa$, and refer to the equivalence classes as
\emph{$\kappa$-classes}.

This paper is organized as follows. After terminology and background
in Section~\ref{sec:prelim}, we show in Section~\ref{sec:bij} how to
associate a poset to each $\kappa$-class, and we characterize
structural properties of these posets. This helps us better understand
admissible sequences as a whole, culminating in a bijection
\begin{equation}
  \label{eq:bijection1}
  \Theta\colon\Acyc(\Gamma)/\simkappa\,\longrightarrow
  \big(\Acyc(\Gamma'_e)/\simkappa\!\big)\,\bigcup\,
  \big(\Acyc(\Gamma''_e)/\simkappa\!\big)\,,
\end{equation}
where $\Gamma'_e$ and $\Gamma''_e$ are the graphs formed by deleting
and contracting a cycle-edge $e$ of $\Gamma$, respectively. From this
bijection, the recursion relation for $\kappa(\Gamma)$
in~\cite{Macauley:08b} becomes an immediate corollary, enumerating
$\kappa(\Gamma)$ through an evaluation of the Tutte polynomial. In
Section~\ref{sec:invariant}, we use our bijection to construct a
complete invariant of $\Acyc(\Gamma)/\simkappa$, the set of
$\kappa$-classes of $\Gamma$. In Section~\ref{sec:conjugacy}, we
review a connection to Coxeter theory, and show how the prior results
easily solve the conjugacy problem for Coxeter elements in all
simply-laced Coxeter groups, and how $\kappa(\Gamma)$ enumerates the
conjugacy classes of Coxeter elements. Finally, in the summary, we
briefly discuss how the equivalence relation~$\tsimkappa$ arises in
other areas of mathematics such as sequential dynamical systems, the
chip-firing game, and the representation theory of quivers. Throughout
the paper, we maintain a running example (that we visit five times)
using a six-vertex graph $\Gamma$ that should enhance the paper's
readability and motivate the main ideas.


\section{Terminology and Background.}
\label{sec:prelim}

Let $\Gamma$ be an undirected, simple and loop-free graph with vertex
set $\vset[\Gamma]=\{1,2,\dots,n\}$ and edge set $\eset[\Gamma]$. Let
$S_\Gamma$ denote the set of total orders (\ie, permutations) of
$\vset[\Gamma]$. Define a relation $\sim$ on $S_\Gamma$ where
$\pi\sim\pi'$ if $\pi=\pi_1\pi_2\cdots\pi_n$ and
$\pi'=\pi'_1\pi'_2\cdots\pi'_n$ differ by a single adjacent
transposition $\pi_i\pi_{i+1}\mapsto\pi_{i+1}\pi_i$ where
$\{\pi_i,\pi_{i+1}\}\not\in\eset[\Gamma]$. The reflexive transitive
closure of $\sim$ is an equivalence relation on $S_\Gamma$ denoted by
$\tsimalpha$. We denote the equivalence class containing $\pi$ by
$[\pi]_\Gamma$, and set
\begin{equation*}
  \SYmod\, = \big\{[\pi]_\Gamma\mid\pi\in S_\Gamma\big\}\,.
\end{equation*}
This corresponds to partially commutative monoids as defined
in~\cite{Cartier:69}, but restricted to fixed length permutations over
$\vset[\Gamma]$ and with commutation relations encoded by
non-adjacency in the graph $\Gamma$. Those familiar with Coxeter
theory will recognize the similarity of these equivalence classes and
the \emph{commutation classes}~\cite{Stembridge:96} of reduced
expressions of Coxeter elements.

Orientations of $\Gamma$ are represented as maps
$O_\Gamma\colon\eset[\Gamma] \longrightarrow \vset[\Gamma] \times
\vset[\Gamma]$, which may also be viewed as directed graphs. The set
of acyclic orientations of $\Gamma$ is denoted by $\Acyc(\Gamma)$, and we
set $\acyc(\Gamma) = \card{\Acyc(\Gamma)}$.
Each acyclic orientation defines a partial ordering on $\vset[\Gamma]$
where the covering relations are $i\leq_{O_\Gamma}\!j$ if $\{i,j\}\in
\eset[\Gamma]$ and $O_\Gamma(\{i,j\}) = (i,j)$.
The set of linear extensions of $O_\Gamma$ contains precisely the
permutations $\pi\in S_\Gamma$ such that if $i\leq_{O_\Gamma} j$, then
$i$ precedes $j$ in $\pi$. Through the ordering of $\vset[\Gamma]$,
every permutation $\pi\in S_\Gamma$ induces a canonical linear order
on $\vset[\Gamma]$. Moreover, each permutation $\pi\in S_\Gamma$
induces an acyclic orientation $O_\Gamma^\pi\in\Acyc(\Gamma)$ defined
by $O_\Gamma^{\pi}(\{i,j\}) = (i,j)$ if $i$ precedes $j$ in $\pi$ and
$O_\Gamma^{\pi}(\{i,j\}) = (j,i)$ otherwise. The canonical bijection
\begin{equation}
  \label{eq:bij1}
  f_\Gamma\colon\SYmod\,\longrightarrow \Acyc(\Gamma)\,,\qquad
     f_\Gamma([\pi]_\Gamma) = O_\Gamma^\pi\,,
\end{equation}
identifies equivalence classes and acyclic orientations, and thus the
number of equivalence classes under $\tsimalpha$ is $\acyc(\Gamma)$.

For $O_\Gamma\in\Acyc(\Gamma)$ and $e = \{v,w\}\in\eset[\Gamma]$, let
$O^{\rho(e)}_\Gamma$ be the orientation of $\Gamma$ obtained from
$O_\Gamma$ by reversing the orientation of the edge $e$. Let
$\Gamma'_e$ and $\Gamma''_e$ denote the graphs obtained from $\Gamma$
by deletion and contraction (see, e.g.,~\cite[p.~415]{Miller:07}) of
$e$, respectively, and let $O_{\Gamma'}$ and $O_{\Gamma''}$ denote the
orientations of $O_\Gamma$ inherited under these operations. (Since
our graphs are assumed to be loop-free, when we contract an edge
$\{v,w\}$, we remove the resulting loop.) The bijection
\begin{equation}
  \label{eq:bijection2}
  \beta_e \colon \Acyc(\Gamma)\longrightarrow\Acyc(\Gamma_e')\cup\Acyc(\Gamma''_e)
\end{equation}
defined by
\begin{eqnarray}
  \label{eq:beta}
  O_\Gamma&\stackrel{\beta_e}{\longmapsto}&
  \begin{cases}
    O_{\Gamma'}\,,  & O^{\rho(e)}_\Gamma\not\in\Acyc(\Gamma)\,,  \\
    O_{\Gamma'}\,,  & O^{\rho(e)}_\Gamma\in\Acyc(\Gamma)\mbox{ and } O_\Gamma(e)=(v,w) \,,  \\
    O_{\Gamma''}\,, & O^{\rho(e)}_\Gamma\in\Acyc(\Gamma)\mbox{ and } O_\Gamma(e)=(w,v) \,,  \\
  \end{cases}
\end{eqnarray}
is well-known, and shows that one may compute $\alpha(\Gamma)$ through the
recursion relation
\begin{equation*}
  \label{eq:alpha}
  \acyc(\Gamma)=\acyc(\Gamma'_e)+\acyc(\Gamma''_e) \,,
\end{equation*}
valid for any $e\in\eset[\Gamma]$. It basically removes the edge
$e=\{v,w\}$ if it cannot be contracted, and otherwise, it either
contracts or removes it depending on its orientation. We illustrate
this with the following example, which we will revisit four more times
throughout this article.

\begin{ex}
\label{ex:beta}
For an explicit example of $\beta_e$, see Figure~\ref{fig:beta}, which
shows three acyclic orientations of the same graph $\Gamma$, and a
fixed edge $e=\{v,w\}$. (The vertices $v'$, $w'$, and $z$, which will
be referred to later, are only labeled once for clarity.) Call these
orientations $O^a_\Gamma$, $O^b_\Gamma$, and $O^c_\Gamma$,
respectively. The map $\beta_e$ removes edge $e$ from $O^a_\Gamma$
because contracting it would result in a directed cycle. Neither
$O^b_\Gamma$ nor $O^c_\Gamma$ have a directed path from $v$ to $w$
other than the edge $(v,w)$, so in both cases, contracting $e$ would
give an acyclic orientation of $\Gamma''_e$. By the definition of
$\beta_e$ in~\eqref{eq:beta}, $\beta_e$ removes $e$ in $O^b_\Gamma$,
and contracts it in $O^c_\Gamma$.
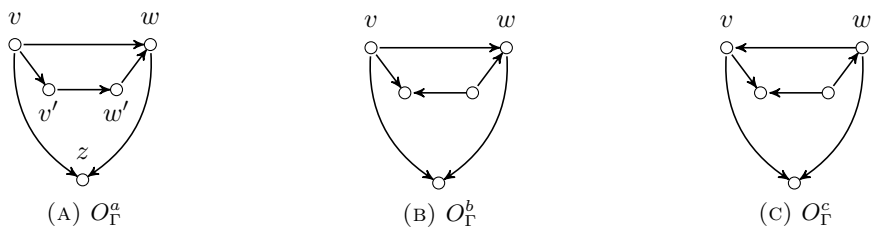
\begin{figure}
\centering
\subfloat[$O_\Gamma^a$] {
\begin{tikzpicture}[scale=0.6]
  \tikzstyle{VertexStyle}=[shape=circle,minimum size=2.5pt,inner sep=1.6pt,draw]
  \Vertex[style={fill=white}, x=1.5, y=4, L=\tiny {}]{v}
  \Vertex[style={fill=white}, x=4.5, y=4, L=\tiny {}]{w}
  \Vertex[style={fill=white}, x=2.25, y=3, L=\tiny {}]{vv}
  \Vertex[style={fill=white}, x=3.75, y=3, L=\tiny {}]{ww}
  \Vertex[style={fill=white}, x=3, y=1, L=\tiny {}]{z}
  \tikzstyle{LabelStyle}=[above]
  \Edge[style = {post}](v)(w)
  \Edge[style = {post}](v)(vv)
  \Edge[style = {post}](vv)(ww)
  \Edge[style = {post}](ww)(w)
  \Edge[style = {post, bend right}](v)(z)
  \Edge[style = {post, bend left}](w)(z)
  \draw[fill=white] (1.5,4) node[label=above:\small $v$]{};
  \draw[fill=white] (4.5,4) node[label=above:\small $w$]{};
  \draw[fill=white] (2.25,3.2) node[label=below:\small $v'$]{};
  \draw[fill=white] (3.75,3.2) node[label=below:\small $w'$]{};
  \draw[fill=white] (3,1) node[label=above:\small $z$]{};
\end{tikzpicture}
}
\hspace{0.8in}
\subfloat[$O_\Gamma^b$] {
\begin{tikzpicture}[scale=0.6]
  \tikzstyle{VertexStyle}=[shape=circle,minimum size=2.5pt,inner sep=1.6pt,draw]
  \Vertex[style={fill=white}, x=1.5, y=4, L=\tiny {}]{v}
  \Vertex[style={fill=white}, x=4.5, y=4, L=\tiny {}]{w}
  \Vertex[style={fill=white}, x=2.25, y=3, L=\tiny {}]{vv}
  \Vertex[style={fill=white}, x=3.75, y=3, L=\tiny {}]{ww}
  \Vertex[style={fill=white}, x=3, y=1, L=\tiny {}]{z}
  \tikzstyle{LabelStyle}=[above]
  \Edge[style = {post}](v)(w)
  \Edge[style = {post}](v)(vv)
  \Edge[style = {pre}](vv)(ww)
  \Edge[style = {post}](ww)(w)
  \Edge[style = {post, bend right}](v)(z)
  \Edge[style = {post, bend left}](w)(z)
  \draw[fill=white] (1.5,4) node[label=above:\small $v$]{};
  \draw[fill=white] (4.5,4) node[label=above:\small $w$]{};
\end{tikzpicture}
}
\hspace{0.8in}
\subfloat[$O_\Gamma^c$] {
\begin{tikzpicture}[scale=0.6]
  \tikzstyle{VertexStyle}=[shape=circle,minimum size=2.5pt,inner sep=1.6pt,draw]
  \Vertex[style={fill=white}, x=1.5, y=4, L=\tiny {}]{v}
  \Vertex[style={fill=white}, x=4.5, y=4, L=\tiny {}]{w}
  \Vertex[style={fill=white}, x=2.25, y=3, L=\tiny {}]{vv}
  \Vertex[style={fill=white}, x=3.75, y=3, L=\tiny {}]{ww}
  \Vertex[style={fill=white}, x=3, y=1, L=\tiny {}]{z}
  \tikzstyle{LabelStyle}=[above]
  \Edge[style = {pre}](v)(w)
  \Edge[style = {post}](v)(vv)
  \Edge[style = {pre}](vv)(ww)
  \Edge[style = {post}](ww)(w)
  \Edge[style = {post, bend right}](v)(z)
  \Edge[style = {post, bend left}](w)(z)
  \draw[fill=white] (1.5,4) node[label=above:\small $v$]{};
  \draw[fill=white] (4.5,4) node[label=above:\small $w$]{};
\end{tikzpicture}
}
\caption{An example of the map $\beta_e$ applied to three acyclic
  orientations of a graph $\Gamma$. If contracting the edge
  $e=\{v,w\}$ would introduce a directed cycle (as in $O_\Gamma^a$),
  then we must delete it. Otherwise, we can either delete or contract
  it, so we pick the convention that we delete it if it is oriented
  $(v,w)$ (as in $O_\Gamma^b$), and contract it if is oriented $(w,v)$
  (as in $O_\Gamma^c$).}
\label{fig:beta}
\end{figure}
\end{ex}

Via the bijection in~\eqref{eq:bij1}, it is clear that mapping
$\pi=\pi_1\pi_2\cdots\pi_n\in[\pi]_\Gamma$ to $\pi_2\cdots\pi_n\pi_1$
corresponds precisely to converting the source vertex $\pi_1$ in
$O_\Gamma^\pi$ into a sink. We call such a conversion a
\emph{source-to-sink operation}, or a \emph{click}.
Two orientations $O_\Gamma,O_\Gamma'\in \Acyc(\Gamma)$ where
$O_\Gamma$ can be transformed into $O'_\Gamma$ by a sequence of clicks
are said to be click-related. We write this as
$\mathbf{c}(O_\Gamma)=O'_\Gamma$ where $\mathbf{c}=c_1c_2\cdots c_k$
with $c_i\in\vset[\Gamma]$. To clarify notation, we mean
\[
\mathbf{c}(O_\Gamma)=c_k(c_{k-1}(\cdots c_2(c_1(O_\Gamma))))\,.
\]
Such as sequence $\mathbf{c}$ is called an \emph{admissible sequence},
or a \emph{click-sequence}. The former term comes from the
representation theory of quivers~\cite{Auslander:97, Kleiner:07}, but
we will usually stick to the latter due to brevity, the overuse of the
term ``admissible sequence'' throughout mathematics, and the
convenience of ``click'' doubling as a verb. It is straightforward to
verify that this click-relation is an equivalence relation on
$\Acyc(\Gamma)$, and we also refer to click-related acyclic
orientations as \emph{$\kappa$-equivalent}. Clearly, and as pointed
out by V.~Reiner~\cite[p. 309]{Novik:02}, one may also approach this
in the setting of total orders on $\vset[\Gamma]$ by identifying
elements that differ by $(i)$ flips of adjacent elements not connected
in $\Gamma$ and $(ii)$ cyclic shifts. However, for our purposes,
approaching this at the level of acyclic orientations seems more
natural in light of the bijection~\eqref{eq:bij1}.


\section{Constructing the Bijection $\Theta$.}
\label{sec:bij}

The bijection $\beta_e \colon \Acyc(\Gamma) \longrightarrow
\Acyc(\Gamma'_e) \cup \Acyc(\Gamma''_e)$ in~\eqref{eq:bijection2} does
not extend to a well-defined map on $\kappa$-classes, i.e.,
$\Acyc(\Gamma)/\simkappa\longrightarrow
(\Acyc(\Gamma'_e)/\simkappa)\cup(\Acyc(\Gamma''_e)/\simkappa)$. Thus, we
need to take a different approach to construct our bijection. An edge
$e$ of an undirected graph $\Gamma$ is a \emph{bridge} if removing $e$
increases the number of connected components of $\Gamma$. An edge that
is not a bridge is a \emph{cycle-edge}, or equivalently, an edge $e$
is a cycle-edge if it is contained in a cycle traversing $e$ precisely
once.

Throughout, we will let $e=\{v,w\}$ denote a fixed cycle-edge of the
connected graph $\Gamma$, and, for ease of notation, we set
$\Gamma'=\Gamma'_e$ and $\Gamma''=\Gamma''_e$. Recall that for
$O_\Gamma\in\Acyc(\Gamma)$ we let $O_{\Gamma'}$ and $O_{\Gamma''}$
denote the inherited orientations of $\Gamma'$ and $\Gamma''$. Notice
that $O_{\Gamma'}$ is always acyclic, while $O_{\Gamma''}$ is acyclic
if and only if there is no directed path with endpoints $v$ and $w$ in
$O_{\Gamma'}$. Finally, we let $[O_\Gamma]$ denote the $\kappa$-class
containing $O_\Gamma$.

The interval $[a,b]$ of a poset $\mathcal{P}$ (where $a\leq b$) is the
subposet consisting of all $c\in \mathcal{P}$ such that $a\leq c\leq
b$. Viewing a finite poset $\mathcal{P}$ as a directed graph
$D_{\mathcal{P}}$, the interval $[a,b]$ contains precisely the
vertices that lie on a directed path from $a$ to $b$, and thus is a
vertex-induced subgraph of $D_{\mathcal{P}}$. By assumption,
$\eset[\Gamma]$ contains $\{v,w\}$, so for all
$O_\Gamma\in\Acyc(\Gamma)$ either $v\leq_{O_\Gamma}\!w$ or
$w\leq_{O_\Gamma}\!v$. In this section, we will study the interval
$[v,w]$ in the poset $O_\Gamma$ (when $v\leq_{O_\Gamma} w$) and its
behavior under clicks.

\begin{defn}
  Let $\Acyc_\le(\Gamma)$ be the set of acyclic orientations of
  vertex-induced subgraphs of $\Gamma$. Define the map
  \begin{equation*}
    \mathcal{I}\colon\Acyc(\Gamma)\longrightarrow\Acyc_\le(\Gamma)
  \end{equation*}
  by $\mathcal{I}(O_\Gamma) = [v,w]$ if $v\leq_{O_\Gamma}\!w$, and by $\mathcal{I}(O_\Gamma)
  = \varnothing$ (the null graph) otherwise. We will refer to
$\mathcal{I}(O_\Gamma)$ as the \emph{$vw$-interval} of $O_\Gamma$.
\end{defn}
Elements of $\Acyc(\Gamma)$ can be thought of as posets over
$\vset[\Gamma]$, and elements of $\Acyc_\le(\Gamma)$ can be thought of
as certain subposets of these, though they need not be induced
(because two vertices on a directed path in $\Gamma$ need not be on a
directed path in an induced subgraph of $\Gamma$). Through a slight
abuse of notation, we will at times refer to $\mathcal{I}(O_\Gamma)$
as a poset, a directed graph, or a subset of $\vset[O_\Gamma]$. In
this last case, it is understood that the relations are inherited from
$O_\Gamma$.

\medskip

Let $P$ be an undirected path in $\Gamma$ of length-$k$, i.e.,
$P=v_0,v_1,\dots,v_{k-1},v_k$ where $\{v_{i-1},v_i\}\in\eset[\Gamma]$
for $i=1,\dots,k$. Define the function
\begin{equation}
  \label{eq:nu_P}
\nu_P\colon\Acyc(\Gamma)\longrightarrow\Z \,,
\end{equation}
where $\nu_P(O_\Gamma)$ is the number of edges in $\Gamma$ of the form
$\{v_{i-1},v_i\}$ oriented as $(v_{i-1},v_i)$ in $O_\Gamma$, minus the
number of edges oriented as $(v_i,v_{i-1})$. If $P$ is a cycle (i.e.,
$v_0=v_k$), $\nu_P$ is preserved under clicks, and thus in this case,
it extends to a map $\nu_P^* \colon \Acyc(\Gamma)/\simkappa
\longrightarrow \Z$. In~\cite{Shi:01}, J.-Y.~Shi defines this function
for Coxeter graphs containing a single cycle, referring to it as
Coleman's \emph{$\nu$-function} (see~\cite{Coleman:89}). The
definition given here is more general, and will allow us to extend
Shi's characterization of conjugacy classes to include all
simply-laced Coxeter groups.

\medskip

\begin{ex}
  \label{ex:nu}
  Continuing with Example~\ref{ex:beta}, consider the three
  orientations in Figure~\ref{fig:beta}, whose $vw$-intervals are the
  following:
  \[
  \mathcal{I}(O^a_\Gamma)=\{v,w,v',w'\}\,,\qquad
  \mathcal{I}(O^b_\Gamma)=\{v,w\}\,,\qquad \mathcal{I}(O^c_\Gamma)=\varnothing\,.
  \]
  Next, consider the undirected (but oriented) path $P=v,v',w',w,v$ (a
  cycle) and the corresponding map $\nu_P$, as defined
  in~\eqref{eq:nu_P}. It is easy to check that
  \[
  \nu_P(O^a_\Gamma)=2\,,\qquad \nu_P(O^b_\Gamma)=0\,,\qquad
  \nu_P(O^c_\Gamma)=2\,.
  \]
  We conclude that $O^b_\Gamma$ cannot be $\kappa$-equivalent to
  either $O^a_\Gamma$ or $O^c_\Gamma$. Finally, if we consider the
  undirected path $Q=v,z,w,v$, we have
  \[
  \nu_Q(O^a_\Gamma)=-1\,,\qquad \nu_Q(O^b_\Gamma)=-1\,,\qquad
  \nu_Q(O^c_\Gamma)=1\,.
  \]
  Therefore, $O^a_\Gamma\nkeq O^c_\Gamma$, and hence all three of
  these orientations lie in distinct $\kappa$-classes.
\end{ex}

As we will see in Section~\ref{sec:invariant}, when taken over all
cycles of $\Gamma$, the $\nu$-function is a actually a complete
invariant, i.e., $O\tsimkappa O'$ if and only iff $\nu_C(O)=\nu_C(O')$
for all cycles $C$ in $\Gamma$. First, we need to establish a series
of structural results about the $vw$-interval. Since
$\{v,w\}\in\eset[\Gamma]$, every $\kappa$-class contains at least one
orientation $O_\Gamma$ with $v\leq_{O_\Gamma}\!w$, and thus there is
at least one element $O_\Gamma$ in each $\kappa$-class with
$\mathcal{I}(O_\Gamma) \ne \varnothing$. As the next result shows,
this (non-empty) choice of $vw$-interval is independent of the choice
of representative from $[O_\Gamma]$, meaning that there is a
well-defined notion of the $vw$-interval of a $\kappa$-equivalence
class. We formalize this by extending the map
$\mathcal{I}\colon\Acyc(\Gamma)\to\Acyc_\leq(\Gamma)$ to a map
$\mathcal{I}^*\colon\Acyc(\Gamma)/\simkappa\to\Acyc_\le(\Gamma)$.
\begin{prop}
  \label{prop:I*}
  The map $\mathcal{I}$ can be extended to a map
  \begin{equation*}
    \mathcal{I}^*\colon\Acyc(\Gamma)/\simkappa\longrightarrow\Acyc_\le(\Gamma)
    \text{\quad by\quad}
    \mathcal{I}^*([O_\Gamma]) = \mathcal{I}(O^1_\Gamma) \,,
  \end{equation*}
  where $O^1_\Gamma$ is any element of $[O_\Gamma]$ for which $\mathcal{I}(O_\Gamma^1)
  \ne \varnothing$.
\end{prop}

\begin{proof}
  It suffices to prove that $\mathcal{I}^*$ is well-defined. Consider
  $O^1_\Gamma\tsimkappa O^2_\Gamma$ with $v\leq_{O^i_\Gamma}\!w$ for
  $i=1,2$. Clearly, $\mathcal{I}(O^1_\Gamma)$ and
  $\mathcal{I}(O^2_\Gamma)$ contain $v$ and $w$, so suppose that
  $a\in\mathcal{I}(O^1_\Gamma)\setminus\{v,w\}$.  Then $a$ lies on a
  directed path $P'$ from $v$ to $w$ in $O^1_\Gamma$, of length $k\geq
  2$ (i.e., $P'$ traverses at least $2$ edges). Let $P$ be the cycle
  formed by adding vertex $v$ to the end of $P'$. Clearly
  $\nu_P(O_\Gamma^1) = k-1$ since $O^1_\Gamma(e) = (v,w)$.

  By assumption, $O^2_\Gamma\in[O^1_\Gamma]$ with
  $v\leq_{O^2_\Gamma}\!w$. Since $\nu_P$ is constant on $[O^1_\Gamma]$
  it follows from $\nu_P(O_\Gamma^1) = k-1 = \nu_P(O_\Gamma^2)$ that
  every edge of $P'$ is oriented identically in $O^1_\Gamma$ and
  $O^2_\Gamma$, and hence that every directed path $P'$ in
  $O_\Gamma^1$ is contained in $O_\Gamma^2$ as well. Therefore,
  $a\in\mathcal{I}(O^2_\Gamma)$, and the reverse inclusion follows by
  an identical argument.
\end{proof}

In light of Proposition~\ref{prop:I*}, we define the $vw$-interval of
a $\kappa$-class $[O_\Gamma]$ to be $\mathcal{I}^*([O_\Gamma])$. The
$vw$-interval will be central in understanding properties of
click-sequences. First, we make a simple observation without proof; it
also appears in~\cite{Speyer:09} in the context of admissible
sequences in Coxeter theory.

\begin{prop}
  \label{prop:alternating}
  Let $O_\Gamma\in\Acyc(\Gamma)$, let $\mathbf{c}=c_1 c_2\cdots c_m$ be an
  associated click-sequence, and consider any directed edge
  $(v_1,v_2)$ in $O_\Gamma$. Then the occurrences of $v_1$ and $v_2$ in
  $\mathbf{c}$ alternate, with $v_1$ appearing first.
\end{prop}

Because $\{v,w\}\in\eset[\Gamma]$, we can say more about the vertices
in $\mathcal{I}(O_\Gamma)$ that appear between successive instances of
in $v$ and $w$ in a click-sequence.

\begin{prop}
  Let $O_\Gamma\in\Acyc(\Gamma)$, and let $\mathbf{c}=c_1 c_2\cdots c_m$ be an
  associated click-sequence that contains every vertex of
  $\mathcal{I}(O_\Gamma)$ at least once and with $c_1=v$. Then every vertex of
  $\mathcal{I}(O_\Gamma)$ appears in $\mathbf{c}$ before any vertex in
  $\mathcal{I}(O_\Gamma)$ appears twice.
\end{prop}

\begin{proof}
  The proof is by contradiction. Assume the statement is false, and
  let $a\in\mathcal{I}(O_\Gamma)$ be the first vertex whose second instance
  in $\mathbf{c}$ occurs before the first instance of some other
  vertex $z \in \mathcal{I}(O_\Gamma)$. If $a\neq v$, then $a$ is not a
  source in $O_\Gamma$, and there exists a directed edge $(a',a)$. By
  Proposition~\ref{prop:alternating}, $a'$ must appear in $\mathbf{c}$
  before the first instance of $a$, but also between the two first
  instances of $a$. This is impossible, because $a$ was chosen to be
  the first vertex appearing twice in $\mathbf{c}$. That only leaves
  $a=v$, and $v$ must appear twice before the first instance of
  $w$. However, this contradicts the statement of
  Proposition~\ref{prop:alternating} because
  $\{v,w\}\in\eset[\Gamma]$. 
\end{proof}

The next result shows that for any click-sequence $\mathbf{c}$ that
contains every element in $\mathcal{I}(O_\Gamma)$ precisely once, we may assume
without loss of generality that the vertices in $\mathcal{I}(O_\Gamma)$
appear consecutively.

\begin{prop}
  \label{prop:forcedclick}
  Let $O_\Gamma\in\Acyc(\Gamma)$ be an acyclic orientation with
  $v\leq_{O_\Gamma}w$.  If $\mathbf{c} = c_1c_2\cdots c_m$ is an
  associated click-sequence containing precisely one instance of $w$,
  and no subsequent instances of vertices from
  $\mathcal{I}(O_\Gamma)$, then there exists a click-sequence
  $\mathbf{c}' = c'_1c'_2\cdots c'_m$ such that $(i)$ there exists an
  interval $[p,q]$ of $\N$ with $c'_j\in\mathcal{I}(O_\Gamma)$ iff
  $p\leq j\leq q$, and $(ii)$ $\mathbf{c}(O_\Gamma) =
  \mathbf{c}'(O_\Gamma)$.
\end{prop}

\begin{proof}
  We prove the proposition by constructing a desired click-sequence
  $\mathbf{c}''$ from $\mathbf{c}$ through a series of transpositions
  where each intermediate click-sequence $\mathbf{c}'$ satisfies
  $\mathbf{c}(O_\Gamma) = \mathbf{c}'(O_\Gamma)$. Such transpositions
  are said to have property $T$.

  Let $I=\mathcal{I}(O_\Gamma)$, and let $A$ be the set of vertices in $I^c
  = \vset[\Gamma]\setminus I$ that lie on a directed path in $O_\Gamma$ to a
  vertex in $I$ (vertices \emph{above} $I$), and let $B$ be the set of
  vertices that lie on a directed path in $O_\Gamma$ from a vertex in $I$
  (vertices \emph{below} $I$). Let $C$ be the complement of $I\cup
  A\cup B$. Two vertices $c_i,c_j\in A\cup B$ with $i<j$ for which
  there is no element $c_k \in A\cup B$ with $i<k<j$ are said to be
  \emph{tight}. We will investigate when transpositions of tight
  vertices in a click-sequence $\mathbf{c}$ of $O_\Gamma$ has property $T$,
  and we will see that this is always the case if $c_i\in B$ and
  $c_j\in A$.
  Consider the intermediate acyclic orientation after applying
  successive clicks $c_1c_2\cdots c_{i-1}$ to $O_\Gamma$. Obviously, $c_i$
  is a source. At this point, if $c_j$ were not a source, then there
  would be an adjacent vertex $a\in A$ with the edge $\{a,c_j\}$
  oriented $(a,c_j)$. For $c_j$ to be clicked as usual (i.e., as a
  source), $a$ must be clicked first, but this would break the
  assumption that $c_i$ and $c_j$ are tight. Therefore, $c_i$ and
  $c_j$ are both sources at this intermediate step, and so the
  vertices $c_i,c_{i+1},\dots,c_j$ are an independent set of sources,
  and may be permuted in any manner without changing the image of the
  click sequence. Therefore, the transposition of $c_i$ and $c_j$ in
  $\mathbf{c}$ has property $T$, as claimed. By iteratively
  transposing tight pairs in $\mathbf{c}$, we can construct a
  click-sequence with the property that every vertex in $A$ comes
  before every vertex in $B$. In light of this, we may assume without
  loss of generality that $\mathbf{c}$ has this property.

  The next step is to show that we can move all vertices in $A$ before
  $v$, and all vertices in $B$ after $w$ via transpositions having
  property $T$.  Let $a$ be the first vertex in $A$ appearing after
  $v$ in the click sequence $\mathbf{c}$.  We claim that the
  transposition moving $a$ to the position directly preceding $v$ has
  property $T$. This is immediate from the observation that when $v$
  is to be clicked, $a$ is a source as well, by the definition of $A$,
  thus it may be clicked before $v$, without preventing subsequent
  clicks of vertices up until the original position of $a$. Therefore,
  we may one-by-one move the vertices in $A$ that are between $v$ and
  $w$, in front of $v$. An analogous argument shows that we may move
  the vertices in $B$ that appear before $w$ to a position directly
  following $w$. In the resulting click-sequence $\mathbf{c}'$, the
  only vertices between $v$ and $w$ are either in $I$ or $C$. The
  subgraph of the directed graph $O_\Gamma$ induced by $C$ is a disjoint
  union of weakly connected components, and none of the vertices are
  adjacent to $I$. By definition of $A$ and $B$, there cannot exist a
  directed edge $(c,a)$ or $(b,c)$, where $a\in A$, $b\in B$, and
  $c\in C$. Thus for each weakly connected component of $C$, the
  vertices in the component can be moved within $\mathbf{c}'$,
  preserving their relative order, to a position either $(i)$ directly
  after the vertices in $A$ and before $v$, or $(ii)$ directly after
  $w$ and before the vertices of $B$. Call this resulting
  click-sequence $\mathbf{c}''$. As we just argued, all the
  transpositions occurring in the rearrangement
  $\mathbf{c}\mapsto\mathbf{c}''$ has property $T$, and
  $\mathbf{c}''$ contains all of the vertices in $I$ in consecutive
  order, and this proves the result. 
\end{proof}

We remark that the last two results together imply that for the
interval $[p,q]$ in the statement of
Proposition~\ref{prop:forcedclick}, $c_p=v$, $c_q=w$, and the sequence
$c_pc_{p+1}\cdots c_q$ contains every vertex in $\mathcal{I}(O_\Gamma)$
precisely once. A simple induction argument implies the following.

\begin{cor}
  \label{cor:intervals}
  Suppose that $O_\Gamma\in\Acyc(\Gamma)$ with $v\leq_{O_\Gamma} w$,
  and let $\mathbf{c}=c_1c_2\cdots c_m$ be a click-sequence where $w$
  appears exactly $k$ times, and no vertex from
  $\mathcal{I}(O_\Gamma)$ appears in $\mathbf{c}$ after the last
  instance of $w$. Then there exists a click-sequence
  $\mathbf{c}'=c'_1c'_2\cdots c'_m$ such that $(i)$ there are $k$
  disjoint intervals $[p_i,q_i]$ of $\N$ such that
  $c_j\in\mathcal{I}(O_\Gamma)$ iff $p_i\leq j\leq q_i$ for some $i$,
  and $(ii)$ $\mathbf{c}(O_\Gamma)=\mathbf{c}'(O_\Gamma)$.
\end{cor}

\begin{proof}
  The argument is by induction on $k$. When $k=1$, the statement is
  simply Proposition~\ref{prop:forcedclick}. Suppose the statement holds
  for all $k\leq N$, for some $N\in\N$, and let $\mathbf{c}$ be a
  click-sequence containing $N+1$ instances of $w$. Let $c_\ell$ be
  the second instance of $v$ in $\mathbf{c}$, and consider the two
  click-sequences $\mathbf{c}_i:=c_1c_2\cdots c_{\ell-1}$ and
  $\mathbf{c}_f:=c_\ell c_{\ell+1}\cdots c_m$. By
  Proposition~\ref{prop:forcedclick}, there exists an interval
  $[p_1,q_1]$ with $p_1<q_1<\ell$, and by the induction hypothesis, there
  exists $k$ intervals $[p_2,q_2],\dots,[p_{k+1},q_{k+1}]$ with
  $\ell\leq p_2<q_2<\cdots<p_{k+1}<q_{k+1}$ such that if
  $c_j\in\mathcal{I}(O_\Gamma)$, then $p_i\leq j\leq q_i$ for some
  $i=1,\dots,k+1$. 
\end{proof}

Let $\varepsilon\colon\Acyc(\Gamma)\longrightarrow\Acyc(\Gamma')$ be
the canonical map that sends $O_\Gamma$ to $O_{\Gamma'}$. This extends
naturally to a map $\varepsilon^*\colon \Acyc(\Gamma)/\simkappa
\longrightarrow \Acyc(\Gamma')/\simkappa$ between $\kappa$-classes.
Define
\begin{equation*}
  \ieext \colon \Acyc(\Gamma')/\simkappa \longrightarrow\Acyc_\le(\Gamma)
\end{equation*}
by $\ieext([O_{\Gamma'}]) = \mathcal{I}(O^1_\Gamma)$ for any
$O^1_\Gamma \in [O_\Gamma]$ such that $\varepsilon^*([O_\Gamma]) =
[O_{\Gamma'}]$ with $\card{\mathcal{I}(O^1_\Gamma)}\geq 3$, and
$\ieext([O_{\Gamma'}])=(v,w)$ (that is, the subgraph induced by
$\{v,w\}$) if no such acyclic orientation $O^1_\Gamma$ exists. The
following result relates the $vw$-intervals of
$\Acyc(\Gamma)/\simkappa$ and $\Acyc(\Gamma')/\simkappa$ through a
commutative diagram involving $\mathcal{I}^*$ and $\ieext$. An
explicit example immediately follows the proof.

\begin{prop}
  \label{prop:diagram}
  The map $\ieext$ is well-defined, and the diagram
  \begin{equation*}
    \xymatrix{
      \Acyc(\Gamma)/\simkappa
      \ar@{->}[rr]^{\mathcal{I}^*}\ar@{->}[dd]_{\varepsilon^*}
      && \Acyc_\le(\Gamma) \\ \\
      \Acyc(\Gamma')/\simkappa\ar@{-->}[uurr]_{\ieext}
    }
  \end{equation*}
  commutes.
\end{prop}

\begin{proof}
  Let $[O_{\Gamma'}] \in\Acyc(\Gamma')/\simkappa$. If there is at most
  one orientation $O_{\Gamma}\in\Acyc(\Gamma)$ such that
  $\card{\mathcal{I}(O_\Gamma)}\geq 3$ and $\varepsilon(O_\Gamma) \in
  [O_{\Gamma'}]$, or if all orientations of the form $O_\Gamma^1$ in
  the definition of $\ieext$ are $\kappa$-equivalent, then both
  statements of the proposition are clear. Assume therefore that there
  are acyclic orientations $O_\Gamma^{\pi}, O_\Gamma^{\sigma} \in
  \Acyc(\Gamma)$ with $O_\Gamma^\pi\nkeq O_\Gamma^\sigma$, but
  $\eta^*_e([O_\Gamma^\pi])=\eta^*_e([O_\Gamma^\sigma])$ and
  $\card{\mathcal{I}(O_\Gamma^\pi)},\card{\mathcal{I}(O_\Gamma^\sigma)}\ge
  3$.
  It suffices to prove that in this case,
  \begin{equation}
    \label{eq:I_pi=I_sigma}
    \mathcal{I}(O^\pi_\Gamma)=\mathcal{I}(O^\sigma_\Gamma)\,.
  \end{equation}
  This is equivalent to showing that the set of $vw$-paths (directed
  paths from $v$ to $w$) in $O^\pi_{\Gamma'}$ is the same as the set of
  $vw$-paths in $O^\sigma_{\Gamma'}$. From this it will also follow that
  the diagram commutes. By assumption, both of these orientations
  contain at least one $vw$-path. We will consider separately the
  cases when these orientations share or do not share a common
  $vw$-path.

  \textit{Case 1: $O_{\Gamma'}^\pi$ and $O_{\Gamma'}^\sigma$ share no common
    $vw$-path.} Let $P_1$ be a length-$k_1$ $vw$-path in $O_{\Gamma'}^\pi$,
  and let $P_2$ be a length-$k_2$ $vw$-path in
  $O_{\Gamma'}^\sigma$. Suppose that in $O_{\Gamma'}^\pi$ there are $k_2^+$
  edges along $P_2$ oriented from $v$ to $w$, and $k_2^-$ edges
  oriented from $w$ to $v$. Likewise, suppose that in $O_{\Gamma'}^\sigma$
  there are $k_1^+$ edges along $P_1$ oriented from $v$ to $w$, and
  $k_1^-$ edges oriented from $w$ to $v$. If $C=P_1P_2^{-1}$ (the
  cycle formed by traversing $P_1$ followed by $P_2$ in reverse), then
  \begin{equation*}
    \nu_C(O_{\Gamma'}^\pi)= k_1^++k_1^-+k_2^--k_2^+,\qquad
    \nu_C(O_{\Gamma'}^\sigma)= k_1^+-k_1^--k_2^--k_2^+\,.
  \end{equation*}
  Equating these values yields $k_1^-+k_2^-=0$, and since these are
  non-negative integers, $k_1^-=k_2^-=0$. We conclude that $P_1$ is a
  $vw$-path in $O_{\Gamma'}^\sigma$ and $P_2$ is a $vw$-path in
  $O_{\Gamma'}^\pi$, contradicting the assumption that $O_{\Gamma'}^\pi$ and
  $O_{\Gamma'}^\sigma$ share no common $vw$-paths.

  \textit{Case 2: $O_{\Gamma'}^\pi$ and $O_{\Gamma'}^\sigma$ share a
    common $vw$-path $P_1$, say of length $k_1$}. If these are the
  only $vw$-paths, we are done. Otherwise, assume without loss of
  generality that $P_2$ is another $vw$-path in $O_{\Gamma'}^\pi$, say
  of length $k_2$. Then if $C=P_1P_2^{-1}$, we have
  $\nu_C(O_{\Gamma'}^\pi)=k_1-k_2$, and hence
  $\nu_C(O_{\Gamma'}^\sigma)=k_1-k_2$. Therefore, $P_2$ is a $vw$-path
  in $O_{\Gamma'}^\sigma$ as well. Because $P_2$ was arbitrary, we
  conclude that $O_{\Gamma'}^\pi$ and $O_{\Gamma'}^\sigma$ share the
  same set of $vw$-paths. Since Case 1 is impossible, we have
  established \eqref{eq:I_pi=I_sigma}, and the proof is complete.
\end{proof}

\begin{ex}
  Consider the orientations $O_\Gamma^a,O_\Gamma^b\in\Acyc(\Gamma)$
  from our running example (see Figure~\ref{fig:beta}). Since the
  $vw$-intervals of $O_\Gamma^a$ and $O_\Gamma^b$ are non-empty (see
  Example~\ref{ex:nu}),
  \[
  \mathcal{I}^*([O_\Gamma^a])=\mathcal{I}(O_\Gamma^a)=\{v,w,v',w'\}\,,
  \qquad\qquad
  \mathcal{I}^*([O_\Gamma^b])=\mathcal{I}(O_\Gamma^b)=\{v,w\}\,.
  \]
  The natural map $\eta_e^*$ simply removes the edge $\{v,w\}$, i.e.,
  \[
  \eta_e^*([O_\Gamma^a])=[O_{\Gamma'}^a]\,,\qquad\qquad
  \eta_e^*([O_\Gamma^b])=[O_{\Gamma'}^b]\,.
  \]
  Finally, Proposition~\ref{prop:diagram} guarantees a well-defined
  map $\ieext$ satisfying $\ieext\circ\eta_e^*=\mathcal{I}^*$, and
  thus
  \[
  \ieext([O_{\Gamma'}^a])=\mathcal{I}^*([O_\Gamma^a])=\{v,w,v',w'\}\,,
  \qquad\qquad
  \ieext([O_{\Gamma'}^b])=\mathcal{I}^*([O_\Gamma^b])=\{v,w\}\,.
  \]
  This is shown in Figure~\ref{fig:diagram}, though note that the
  domains are the actual $\kappa$-classes containing the given
  orientations, not the orientations themselves.
\begin{figure}
\begin{tikzpicture}[scale=0.5]
  \tikzstyle{VertexStyle}=[shape=circle,minimum size=2.5pt,inner sep=1.6pt,draw]
  \Vertex[style={fill=white}, x=1.5, y=4, L=\tiny {}]{v}
  \Vertex[style={fill=white}, x=4.5, y=4, L=\tiny {}]{w}
  \Vertex[style={fill=white}, x=2.25, y=3, L=\tiny {}]{vv}
  \Vertex[style={fill=white}, x=3.75, y=3, L=\tiny {}]{ww}
  \Vertex[style={fill=white}, x=3, y=1, L=\tiny {}]{z}
  \Edge[style = {post}](v)(w)
  \Edge[style = {post}](v)(vv)
  \Edge[style = {post}](vv)(ww)
  \Edge[style = {post}](ww)(w)
  \Edge[style = {post, bend right}](v)(z)
  \Edge[style = {post, bend left}](w)(z)
  \draw[fill=white] (1.5,4) node[label=above:\small $v$]{};
  \draw[fill=white] (4.5,4) node[label=above:\small $w$]{};
  \Vertex[style={fill=white}, x=10.5, y=3, L=\tiny {}]{v2}
  \Vertex[style={fill=white}, x=13.5, y=3, L=\tiny {}]{w2}
  \Vertex[style={fill=white}, x=11.25, y=2, L=\tiny {}]{vv2}
  \Vertex[style={fill=white}, x=12.75, y=2, L=\tiny {}]{ww2}
  \tikzstyle{LabelStyle}=[above]
  \Edge[style = {post}](v2)(w2)
  \Edge[style = {post}](v2)(vv2)
  \Edge[style = {post}](vv2)(ww2)
  \Edge[style = {post}](ww2)(w2)
  \draw[fill=white] (10.5,3) node[label=above:\small $v$]{};
  \draw[fill=white] (13.5,3) node[label=above:\small $w$]{};
  \Vertex[style={fill=white}, x=1.5, y=-3, L=\tiny {}]{v}
  \Vertex[style={fill=white}, x=4.5, y=-3, L=\tiny {}]{w}
  \Vertex[style={fill=white}, x=2.25, y=-4, L=\tiny {}]{vv}
  \Vertex[style={fill=white}, x=3.75, y=-4, L=\tiny {}]{ww}
  \Vertex[style={fill=white}, x=3, y=-6, L=\tiny {}]{z}
  \Edge[style = {post}](v)(vv)
  \Edge[style = {post}](vv)(ww)
  \Edge[style = {post}](ww)(w)
  \Edge[style = {post, bend right}](v)(z)
  \Edge[style = {post, bend left}](w)(z)
  \tikzstyle{VertexStyle}=[shape=circle,minimum size=2.5pt,inner sep=0pt]
  \Vertex[style={}, x=5.25, y=2.5, L=\tiny {}]{tail1}
  \Vertex[style={fill=white}, x=9, y=2.5, L=\tiny {}]{head1}
  \tikzstyle{LabelStyle}=[above]
  \Edge[style = {post}, label=$\mathcal{I}^*$](tail1)(head1)
  \Vertex[style={}, x=3, y=0.25, L=\tiny {}]{tail2}
  \Vertex[style={fill=white}, x=3, y=-2.25, L=\tiny {}]{head2}
  \tikzstyle{LabelStyle}=[left]
  \Edge[style = {post}, label=$\eta_e^*$](tail2)(head2)
  \Vertex[style={}, x=5.25, y=-2.25, L=\tiny {}]{tail3}
  \Vertex[style={fill=white}, x=9, y=1, L=\tiny {}]{head3}
  \tikzstyle{LabelStyle}=[below]
  \Edge[style = {post}, label=$\quad\;\ieext$](tail3)(head3)
  \Edge[style = {post}](tail3)(head3)
\end{tikzpicture}
\hspace{0.85in}
\begin{tikzpicture}[scale=0.5]
  \tikzstyle{VertexStyle}=[shape=circle,minimum size=2.5pt,inner sep=1.6pt,draw]
  \Vertex[style={fill=white}, x=1.5, y=4, L=\tiny {}]{v}
  \Vertex[style={fill=white}, x=4.5, y=4, L=\tiny {}]{w}
  \Vertex[style={fill=white}, x=2.25, y=3, L=\tiny {}]{vv}
  \Vertex[style={fill=white}, x=3.75, y=3, L=\tiny {}]{ww}
  \Vertex[style={fill=white}, x=3, y=1, L=\tiny {}]{z}
  \tikzstyle{LabelStyle}=[above]
  \Edge[style = {post}](v)(w)
  \Edge[style = {post}](v)(vv)
  \Edge[style = {post}](ww)(vv)
  \Edge[style = {pre}](w)(ww)
  \Edge[style = {post, bend right}](v)(z)
  \Edge[style = {pre, bend right}](z)(w)
  \draw[fill=white] (1.5,4) node[label=above:\small $v$]{};
  \draw[fill=white] (4.5,4) node[label=above:\small $w$]{};
  \Vertex[style={fill=white}, x=10.5, y=3, L=\tiny {}]{v2}
  \Vertex[style={fill=white}, x=13.5, y=3, L=\tiny {}]{w2}
  \tikzstyle{LabelStyle}=[above]
  \Edge[style = {post}](v2)(w2)
  \draw[fill=white] (10.5,3) node[label=above:\small $v$]{};
  \draw[fill=white] (13.5,3) node[label=above:\small $w$]{};
  \Vertex[style={fill=white}, x=1.5, y=-3, L=\tiny {}]{v}
  \Vertex[style={fill=white}, x=4.5, y=-3, L=\tiny {}]{w}
  \Vertex[style={fill=white}, x=2.25, y=-4, L=\tiny {}]{vv}
  \Vertex[style={fill=white}, x=3.75, y=-4, L=\tiny {}]{ww}
  \Vertex[style={fill=white}, x=3, y=-6, L=\tiny {}]{z}
  \tikzstyle{LabelStyle}=[above]
  \Edge[style = {post}](v)(vv)
  \Edge[style = {post}](ww)(vv)
  \Edge[style = {pre}](w)(ww)
  \Edge[style = {post, bend right}](v)(z)
  \Edge[style = {pre, bend right}](z)(w)
  \tikzstyle{VertexStyle}=[shape=circle,minimum size=2.5pt,inner sep=0pt]
  \Vertex[style={}, x=5.25, y=2.5, L=\tiny {}]{tail1}
  \Vertex[style={fill=white}, x=9, y=2.5, L=\tiny {}]{head1}
  \tikzstyle{LabelStyle}=[above]
  \Edge[style = {post}, label=$\mathcal{I}^*$](tail1)(head1)
  \Vertex[style={}, x=3, y=0.25, L=\tiny {}]{tail2}
  \Vertex[style={fill=white}, x=3, y=-2.25, L=\tiny {}]{head2}
  \tikzstyle{LabelStyle}=[left]
  \Edge[style = {post}, label=$\eta_e^*$](tail2)(head2)
  \Vertex[style={}, x=5.25, y=-2.25, L=\tiny {}]{tail3}
  \Vertex[style={fill=white}, x=9, y=1, L=\tiny {}]{head3}
  \tikzstyle{LabelStyle}=[below]
  \Edge[style = {post}, label=$\quad\;\ieext$](tail3)(head3)
  \Edge[style = {post}](tail3)(head3)
  \end{tikzpicture}
\caption{An explicit example of the commutative diagram relating the
  maps $\eta_e^*$, $\mathcal{I}^*$, and $\ieext$ from
  Proposition~\ref{prop:diagram}. The domain of these maps are
  actually the sets of $\kappa$-equivalence classes, e.g.,
  $[O_\Gamma^a]$ (left) and $[O_\Gamma^b]$ (right), but they are shown
  acting on actual orientations ($O_\Gamma^a$ and $O_\Gamma^b$) for
  clarity.}
\label{fig:diagram}
\end{figure}
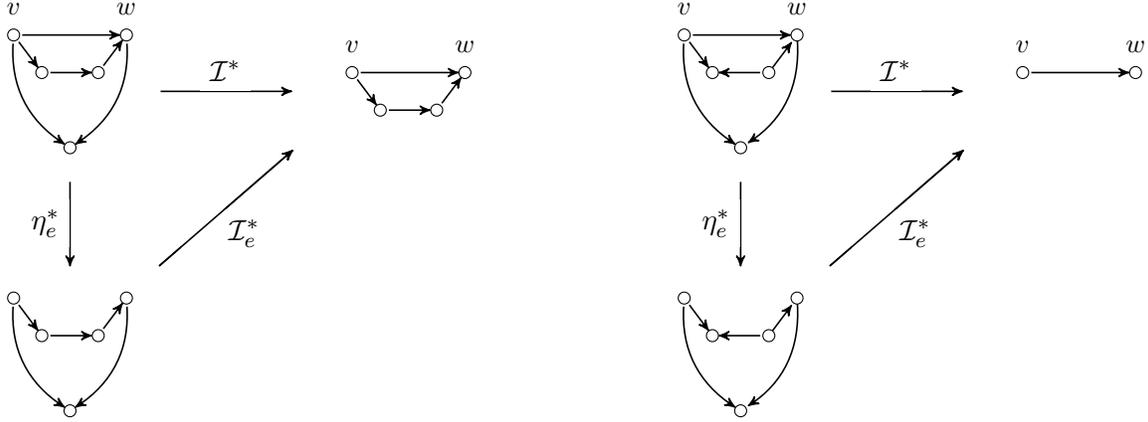
\end{ex}

Let $O_\Gamma\in\Acyc(\Gamma)$ and assume $I = \mathcal{I}(O_\Gamma)$
is non-empty. We write $\Gamma_I$ for the graph formed from $\Gamma$
by contracting all vertices in $I$ to a single vertex, which we denote
by $V_I$. Note that if $I$ only contains $v$ and $w$ then $\Gamma_I =
\Gamma''_e$. Moreover, $O_\Gamma$ gives rise to an orientation
$O_{\Gamma_I}$ of $\Gamma_I$, and this orientation is clearly acyclic.

\begin{prop}
  \label{prop:Y_I}
  Let $O^1_\Gamma, O^2_\Gamma \in \Acyc(\Gamma)$ and assume
  $\mathcal{I}(O^1_\Gamma) = \mathcal{I}(O^2_\Gamma)$. If $O^1_\Gamma
  \nkeq O^2_\Gamma$ then $[O^1_{\Gamma_I}] \nkeq [O^2_{\Gamma_I}]$.
\end{prop}

\begin{proof}
  We prove the contrapositive statement. Set
  $I=\mathcal{I}(O^1_\Gamma)$, suppose $\card{I} = k$, and let
  $v_1v_2\cdots v_k$ be a linear extension of $I$. For any
  click-sequence $\mathbf{c}_I$ between two acyclic orientations
  $O^1_{\Gamma_I}$ and $O^2_{\Gamma_I}$ in $\Acyc(\Gamma_I)$, let
  $\mathbf{c}$ be the click-sequence formed by replacing every
  occurrence of $V_I$ in $\mathbf{c}_I$ by the sequence
  $v_1\cdots v_k$. Then $\mathbf{c}(O^1_\Gamma) = O^2_\Gamma$ and
  $O^1_\Gamma \tsimkappa O^2_\Gamma$ as claimed. 
\end{proof}


We can now utilize the results on poset structure just developed to
establish a bijection
\begin{equation*}
\Theta \colon \Acyc(\Gamma)/\simkappa
\longrightarrow
\bigl(\Acyc(\Gamma'_e)/\simkappa\cup\,\,\Acyc(\Gamma''_e)/\simkappa\bigr)\,,
\end{equation*}
valid for any cycle-edge $e$.
For $[O_\Gamma]\in\Acyc(\Gamma)/\simkappa$, let $O^\pi_\Gamma$ denote
an orientation in $[O_\Gamma]$ such that $\pi=v\pi_2\cdots\pi_n$ and
$w=\pi_i$ for $i$ minimal. We define $\Theta$ by
\begin{equation}
  [O_\Gamma]\stackrel{\Theta}{\longmapsto}
\begin{cases}
  [O^\pi_{\Gamma''}], &
  \mbox{$\exists O^{\pi}_\Gamma\in[O_\Gamma]$ with $\pi=vw\pi_3\cdots\pi_n$} \\
  [O^\pi_{\Gamma'}], & \mbox{otherwise.}
  \end{cases}
\end{equation}

Note that $[O_\Gamma]$ is mapped into $\Acyc(\Gamma'')/\simkappa$ if
and only if the only vertices in $\ieext([O_\Gamma])$ are $v$ and
$w$. Since $\kappa$-equivalence over $\Gamma$ implies
$\kappa$-equivalence over $\Gamma'$, $\Theta$ does not depend on the
choice of $\pi$, and thus is well-defined. We continue our running
example below to illustrate this.
\begin{ex}
\label{ex:Theta}
The orientations $O^a_\Gamma$ and $O^b_\Gamma$ from our running
example are shown in Figure~\ref{fig:Theta}. Since
$\ieext([O^a_\Gamma])=\{v,v',w',w\}$, the map $\Theta$ removes the
edge $e=\{v,w\}$. However, since $\ieext([O^b_\Gamma])=\{v,w\}$, the
map $\Theta$ contracts $e$.
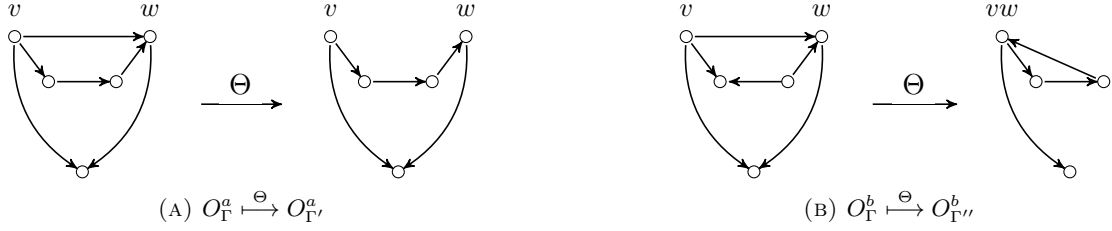
\begin{figure}
\subfloat[$O_\Gamma^a\stackrel{\Theta}{\longmapsto}O_{\Gamma'}^a$] {
\begin{tikzpicture}[scale=0.6]
  \tikzstyle{VertexStyle}=[shape=circle,minimum size=2.5pt,inner sep=1.6pt,draw]
  \Vertex[style={fill=white}, x=1.5, y=4, L=\tiny {}]{v}
  \Vertex[style={fill=white}, x=4.5, y=4, L=\tiny {}]{w}
  \Vertex[style={fill=white}, x=2.25, y=3, L=\tiny {}]{vv}
  \Vertex[style={fill=white}, x=3.75, y=3, L=\tiny {}]{ww}
  \Vertex[style={fill=white}, x=3, y=1, L=\tiny {}]{z}
  \Edge[style = {post}](v)(w)
  \Edge[style = {post}](v)(vv)
  \Edge[style = {post}](vv)(ww)
  \Edge[style = {post}](ww)(w)
  \Edge[style = {post, bend right}](v)(z)
  \Edge[style = {post, bend left}](w)(z)
  \draw[fill=white] (1.5,4) node[label=above:\small $v$]{};
  \draw[fill=white] (4.5,4) node[label=above:\small $w$]{};
  \tikzstyle{VertexStyle}=[shape=circle,minimum size=2.5pt,inner sep=0pt]
  \Vertex[style={}, x=5.5, y=2.5, L=\tiny {}]{a}
  \Vertex[style={fill=white}, x=7.5, y=2.5, L=\tiny {}]{b}
  \tikzstyle{LabelStyle}=[above]
  \Edge[style = {post}, label=$\Theta$](a)(b)
  \tikzstyle{VertexStyle}=[shape=circle,minimum size=2.5pt,inner sep=1.6pt,draw]
  \Vertex[style={fill=white}, x=8.5, y=4, L=\tiny {}]{v2}
  \Vertex[style={fill=white}, x=11.5, y=4, L=\tiny {}]{w2}
  \Vertex[style={fill=white}, x=9.25, y=3, L=\tiny {}]{vv2}
  \Vertex[style={fill=white}, x=10.75, y=3, L=\tiny {}]{ww2}
  \Vertex[style={fill=white}, x=10, y=1, L=\tiny {}]{z2}
  \Edge[style = {post}](v2)(vv2)
  \Edge[style = {post}](vv2)(ww2)
  \Edge[style = {post}](ww2)(w2)
  \Edge[style = {post, bend right}](v2)(z2)
  \Edge[style = {post, bend left}](w2)(z2)
  \draw[fill=white] (8.5,4) node[label=above:\small $v$]{};
  \draw[fill=white] (11.5,4) node[label=above:\small $w$]{};
\end{tikzpicture}
}
\hspace{0.8in}
\subfloat[$O_\Gamma^b\stackrel{\Theta}{\longmapsto}O_{\Gamma''}^b$] {
\begin{tikzpicture}[scale=0.6]
  \tikzstyle{VertexStyle}=[shape=circle,minimum size=2.5pt,inner sep=1.6pt,draw]
  \Vertex[style={fill=white}, x=1.5, y=4, L=\tiny {}]{v}
  \Vertex[style={fill=white}, x=4.5, y=4, L=\tiny {}]{w}
  \Vertex[style={fill=white}, x=2.25, y=3, L=\tiny {}]{vv}
  \Vertex[style={fill=white}, x=3.75, y=3, L=\tiny {}]{ww}
  \Vertex[style={fill=white}, x=3, y=1, L=\tiny {}]{z}
  \Edge[style = {post}](v)(w)
  \Edge[style = {post}](v)(vv)
  \Edge[style = {pre}](vv)(ww)
  \Edge[style = {post}](ww)(w)
  \Edge[style = {post, bend right}](v)(z)
  \Edge[style = {post, bend left}](w)(z)
  \draw[fill=white] (1.5,4) node[label=above:\small $v$]{};
  \draw[fill=white] (4.5,4) node[label=above:\small $w$]{};
  \tikzstyle{VertexStyle}=[shape=circle,minimum size=2.5pt,inner sep=0pt]
  \Vertex[style={}, x=5.5, y=2.5, L=\tiny {}]{a}
  \Vertex[style={fill=white}, x=7.58, y=2.5, L=\tiny {}]{b}
  \tikzstyle{LabelStyle}=[above]
  \Edge[style = {post}, label=$\Theta$](a)(b)
  \tikzstyle{VertexStyle}=[shape=circle,minimum size=2.5pt,inner sep=1.6pt,draw]
  \Vertex[style={fill=white}, x=8.5, y=4, L=\tiny {}]{VW}
  \Vertex[style={fill=white}, x=9.25, y=3, L=\tiny {}]{vv2}
  \Vertex[style={fill=white}, x=10.75, y=3, L=\tiny {}]{ww2}
  \Vertex[style={fill=white}, x=10, y=1, L=\tiny {}]{z2}
  \Edge[style = {post}](VW)(vv2)
  \Edge[style = {post}](vv2)(ww2)
  \Edge[style = {post}](ww2)(VW)
  \Edge[style = {post, bend right}](VW)(z2)
  \draw[fill=white] (8.5,4) node[label=above:\small $vw$]{};
\end{tikzpicture}
}
\caption{An example of the map $\Theta$ applied to the orientations
  $O_\Gamma^a$ and $O_\Gamma^b$ from Example~\ref{ex:beta} and
  Figure~\ref{fig:beta}. If contracting the edge $e=\{v,w\}$ would
  introduce a directed cycle (as in $O_\Gamma^a$), then we must delete
  it. Otherwise (as in $O_\Gamma^b$), contract it. Note that this
  happens precisely when $(v,w)$ is the only directed path from $v$ to
  $w$.  }
\label{fig:Theta}
\end{figure}
\end{ex}

The results we have derived for the $vw$-interval now allow us to
establish the following:
\begin{thm}
  \label{thm:bijection}
  The map $\Theta$ is a bijection.
\end{thm}
\begin{proof}
  We first prove that $\Theta$ is surjective. Let $I=\{v,w\}$ and
  consider an element $[O_{\Gamma''}]\in\Acyc(\Gamma'')/\tsimkappa$
  with $O^{\pi}_{\Gamma''} \in [O_{\Gamma''}]$ where $\pi =
  V_I\pi_2\cdots\pi_{n-1}$. Let $\pi^+=vw\pi_2\cdots\pi_{n-1} \in
  S_\Gamma$.  Clearly $[O^{\pi^+}_\Gamma] \in
  \Acyc(\Gamma)/\tsimkappa$ is mapped to $[O_{\Gamma''}]$ by
  $\Theta$.

  Next, consider an element
  $[O_{\Gamma'}]\in\Acyc(\Gamma')/\tsimkappa$. If there is no element
  $O^\pi_{\Gamma'}$ of $[O_{\Gamma'}]$ such that $\pi =
  vw\pi_3\cdots\pi_n$, then no elements of $[O_\Gamma]$ are of this
  form either, and by definition $[O_{\Gamma'}]$ has a preimage under
  $\Theta$.
  We are left with the case where $[O_{\Gamma'}]$ contains an element
  $O^\pi_{\Gamma'}$ such that $\pi=vw\pi_3\cdots\pi_n$, and we must show
  that there exists $O^{\pi'}_{\Gamma'}\in[O_{\Gamma'}]$ such that
  $[O^{\pi'}_\Gamma]$ contains no element of the form $O^\sigma_\Gamma$ with
  $\sigma = vw\sigma_3\cdots\sigma_n$.
  Note that if $\sigma = vw\sigma_3\cdots\sigma_n$, then the vertices
  in $\mathcal{I}(O^\sigma_\Gamma)$ are precisely $v$ and $w$. If the
  orientation $O_{\Gamma'}$ had a directed path from $v$ to $w$, then
  the corresponding orientation $O_\Gamma\in\Acyc(\Gamma)$ formed by
  adding the edge $e$ with orientation $(v,w)$ has $vw$-interval of
  size at least 3, so by Proposition~\ref{prop:I*}, the acyclic
  orientation $O_\Gamma$ cannot be $\kappa$-equivalent to any
  orientation $O^\sigma_\Gamma$ such that
  $\sigma=vw\sigma_3\cdots\sigma_n$.

  Thus it remains to consider the case when $[O_{\Gamma'}]$ contains no
  acyclic orientation with a directed path from $v$ to $w$. Pick any
  simple undirected path $P'$ from $v$ to $w$ in $\Gamma'$, which is
  possible since $e$ is a cycle-edge. Choose an orientation in
  $[O_{\Gamma'}]$ for which $\nu_{P'}$ is maximal. Without loss of
  generality we may assume that $O_{\Gamma'}$ is this orientation.
  Let $O_\Gamma\in\Acyc(\Gamma)$ be the orientation that agrees with
  $O_{\Gamma'}$, and with $e$ oriented as $(w,v)$. Since we have
  assumed that there is no directed path from $v$ to $w$ this
  orientation is acyclic. We claim that for any
  $\sigma=vw\sigma_3\cdots\sigma_n$ one has $O^\sigma_\Gamma \not\in
  [O_\Gamma]$. To see this, assume the statement is false. Let $P$ be
  the undirected cycle in $\Gamma$ formed by adding the edge $e$ to
  the path $P'$. Because $e$ is oriented as $(v,w)$ in
  $O^\sigma_\Gamma$ and as $(w,v)$ in $O_\Gamma$, we have
  $\nu_P(O^\sigma_\Gamma) = \nu_{P'}(O^\sigma_{\Gamma'})-1$ and
  $\nu_P(O_\Gamma)=\nu_{P'}(O_{\Gamma'})+1$. If $O_\Gamma$ and
  $O^\sigma_\Gamma$ were $\kappa$-equivalent, then
  \begin{equation*}
    \nu_{P'}(O^\sigma_{\Gamma'})-1=\nu_P(O^\sigma_\Gamma)=\nu_P(O_\Gamma)
    =\nu_{P'}(O_\Gamma)+1\,,
  \end{equation*}
  and thus $\nu_{P'}(O^\sigma_{\Gamma'})=\nu_{P'}(O_\Gamma)+2$. Any
  click sequence mapping $O_\Gamma$ to $O^\sigma_\Gamma$ is a
  click-sequence from $O_{\Gamma'}$ to
  $O^\sigma_{\Gamma'}$. Therefore, $O^\sigma_{\Gamma'} \in
  [O_{\Gamma'}]$, which contradicts the maximality of
  $\nu_{P'}(O_{\Gamma'})$. We therefore conclude that $O^\sigma_\Gamma
  \not\in [O_\Gamma]$, that $\Theta([O_\Gamma]) = [O_{\Gamma'}]$, and
  hence that $\Theta$ is surjective.

\medskip

  We next prove that $\Theta$ is an injection. By
  Proposition~\ref{prop:Y_I} (with $I=\{v,w\}$), $\Theta$ is injective
  when restricted to the preimage of $[O_{\Gamma''}]$ under
  $\Theta$. Thus it suffices to show that every element in
  $\Acyc(\Gamma')/\tsimkappa$ has a unique preimage under $\Theta$.
  By Proposition~\ref{prop:diagram}, every preimage of $[O_{\Gamma'}]$
  must have the same $vw$-interval $I$, containing $k>2$
  vertices. Suppose there were preimages
  $[O_\Gamma^\pi]\neq[O_\Gamma^\sigma]$ of $[O_{\Gamma'}]$.  By
  Proposition~\ref{prop:Y_I}, it follows that $O^\pi_{\Gamma_I}\nkeq
  O^\sigma_{\Gamma_I}$. We will now show that this leads to a
  contradiction.

  Assume that $\mathbf{c}=c_1\cdots c_m$ is a click-sequence from
  $O_{\Gamma'}^\pi$ to $O_{\Gamma'}^\sigma$. If one of $\pi$ or
  $\sigma$ is not $\kappa$-equivalent to a permutation with vertices
  $v$ and $w$ in succession, then their corresponding $\kappa$-classes
  would be unchanged by the removal of edge $e$. In light of this, we
  may assume that $\pi = v\pi_2\dots\pi_{n-1}w$ and $\sigma =
  v\sigma_2 \dots \sigma_{n-1}w$, and thus that $c_1 = v$ and $c_m =
  w$. By Proposition~\ref{prop:forcedclick}, we may assume that the
  vertices in $I$ appear in $\mathbf{c}$ in some number of disjoint
  consecutive ``blocks,'' i.e., subsequences of the form $c_i\cdots
  c_{i+k-1}$.  Replacing each of these blocks with $V_I$ yields a
  click-sequence from $O^\pi_{\Gamma_I}$ to $O^\sigma_{\Gamma_{I}}$,
  contradicting the fact that $O^\pi_{\Gamma_I}\nkeq
  O^\sigma_{\Gamma_I}$. Therefore, no such click sequence $\mathbf{c}$
  exists, and $\Theta$ must be an injection, and the proof is
  complete. 
\end{proof}

The main result in~\cite{Macauley:08b} is a recurrence relation for
$\kappa(\Gamma)$ under edge deletion and edge contraction. This is an
easy consequence of Theorem~\ref{thm:bijection}.
\begin{cor}[\cite{Macauley:08b}]
  \label{cor:kappa}
  Let $\Gamma$ be a finite undirected graph with $e\in\eset[\Gamma]$,
  and let $\Gamma_e'$ be the graph obtained from $\Gamma$ by deleting
  $e$, and let $\Gamma_e''$ be the graph obtained from $\Gamma$ by
  contracting $e$. Then
  \begin{equation}
    \label{eq:kappa}
    \kappa(\Gamma) =
      \left\{
      \begin{array}{cl}
        \kappa(\Gamma_1)\kappa(\Gamma_2), &\;\;
          \mbox{$e$ is a bridge linking $\Gamma_1$ and $\Gamma_2$}\,, \\
        \kappa(\Gamma_e')+\kappa(\Gamma_e''), &\;\;
          \mbox{$e$ is a cycle-edge}\,.
      \end{array}
      \right.
  \end{equation}
\end{cor}
The first part involving a bridge is straightforward, while the second
part is a direct consequence of Theorem~\ref{thm:bijection}.


\section{A Complete Invariant of $\Acyc(\Gamma)/\simkappa$.}
\label{sec:invariant}

Theorem~\ref{thm:bijection} is more than just an alternative proof of
the enumeration of $\kappa(\Gamma)$. We can utilize the explicit
bijection to derive an additional interesting and useful corollary:
When taken over all cycles $C$ in a graph $\Gamma$, $\nu_C$ is a
complete invariant of $\Acyc(\Gamma)/\simkappa$. This result is
originally due to Pretzel~\cite{Pretzel:86}, though the techniques are
much different than the ones here.
\begin{thm}
  \label{thm:nu}
  If $\nu_C(O^1_\Gamma)=\nu_C(O^2_\Gamma)$ for every cycle $C$, then
  $O^1_\Gamma\tsimkappa O^2_\Gamma$.
\end{thm}
\begin{proof}
  Assume the statement is false and let $\Gamma$ be a graph for which
  it fails, minimal with respect to $|\eset[\Gamma]|$.  Fix a
  cycle-edge $e=\{v,w\}$, and for any
  $[O_\Gamma]\in\Acyc(\Gamma)/\simkappa$, call an orientation
  $O_\Gamma^\pi\in[O_\Gamma]$ \textit{distinguished} with respect to
  $e$ if $\pi=\pi_1\pi_2\cdots\pi_n$ such that $(i)$ $\pi_1=v$, and
  $(ii)$ $\pi_k=w$ where $k$ is minimal given that $\pi_1=v$. By
  assumption, there exists $O_\Gamma^1\nkeq O_\Gamma^2$ with
  $\nu_C(O_\Gamma^1)=\nu_C(O_\Gamma^2)$ for every cycle $C$ in
  $\Gamma$. Without loss of generality, we may assume that
  $O_\Gamma^1$ and $O_\Gamma^2$ are distinguished orientations with
  respect to $e = \{v,w\}$. Define a $vw$-path to be a directed path
  from $v$ to $w$ that does not traverse $e$. There are three cases to
  consider:

  \textit{Case 1: Both $O_\Gamma^1$ and $O_\Gamma^2$ contain a
  $vw$-path}. By definition of $\Theta$, both $\Theta(O_\Gamma^1)$ and
  $\Theta(O_\Gamma^2)$ are contained in
  $\Acyc(\Gamma'_e)/\simkappa$. Moreover, because $\Theta$ is a
  bijection, they are distinct elements in
  $\Acyc(\Gamma'_e)/\simkappa$. Every cycle in $\Gamma'_e$ is also a
  cycle in $\Gamma$, and therefore,
  $\nu_{C'}(\Theta(O_\Gamma^1))=\nu_{C'}(\Theta(O_\Gamma^2))$ for
  every cycle $C'$ in $\Gamma'_e$, contradicting the minimality of
  $|\eset[\Gamma]|$.

  \textit{Case 2: Neither $O_\Gamma^1$ nor $O_\Gamma^2$ contain a
    $vw$-path}. Since $O_\Gamma^1$ and $O_\Gamma^2$ are distinguished
    with respect to $e$, $\Theta([O_\Gamma^i]) =
    [O_{\Gamma''}^i]\in\Acyc(\Gamma'')/\simkappa$ for $i=1,2$. Again,
    these two orientations are distinct because $\Theta$ is a
    bijection. Any cycle $C''$ in $\Gamma''_e$ beginning and ending at
    the vertex $V$ (the image of $v$ and $w$ under contraction) can be
    canonically extended to a cycle $C$ of $\Gamma$. Therefore,
    $\nu_{C''}(\Theta(O_{\Gamma''}^1))=\nu_{C''}(\Theta(O_{\Gamma''}^2))$,
    again contradicting the minimality of $|\eset[\Gamma]|$.

  \textit{Case 3: Precisely one of $O_\Gamma^1$ and $O_\Gamma^2$
    contain a $vw$-path}. Without loss of generality, suppose that $P$
    is a length-$k$ $vw$-path in $O_\Gamma^1$, and let $C$ be the
    cycle formed by adding vertex $v$ to the end of $P$. Clearly,
    $\nu_C(O_\Gamma^1)=k-1$, and by assumption,
    $\nu_C(O_\Gamma^2)=k-1$ as well. However, this means that every
    edge in $P$ is oriented from $v$ to $w$ in $O_\Gamma^2$,
    contradicting the assumption that $O_\Gamma^2$ did not contain a
    $vw$-path. Therefore, Case 3 is impossible, and the proof is
    complete. 
\end{proof}

\begin{ex}
  \label{ex:cycle} 
  It is immediate from the recurrence~\eqref{eq:kappa} in
  Corollary~\ref{cor:kappa} that $\kappa(\Gamma)=1$ if and only if
  $\Gamma$ is a forest. Let $\Circle_n$ be a chordless $n$-cycle, and
  $\Line_n$ the line graph on $n$ vertices. Deleting or contracting
  any edge of $\Circle_3$ leaves a tree, and so
  $\kappa(\Circle_3)=2$. For $n>3$, deleting an edge from $\Circle_n$
  leaves $\Line_n$, and contracting an edge leaves $\Circle_{n-1}$. By
  Corollary~\ref{cor:kappa}, if $n>3$, then
  \begin{equation}
    \label{eq:Circle_n}
    \kappa(\Circle_n)=\kappa(\Line_n)+\kappa(\Circle_{n-1})
    =1+\kappa(\Circle_{n-1})\,.
  \end{equation}
  From the base case of $\kappa(\Circle_3)=2$, we immediately deduce that
  $\kappa(\Circle_n)=n-1$.
\end{ex}

\begin{ex}
  Let $\Gamma$ be the undirected version of the orientations from our
  running example, with edge $e=\{v,w\}$ as before. Deleting $e$
  leaves $\Gamma'_e$, a $5$-cycle. Contracting $e$ leaves
  $\Gamma''_e$, a $3$-cycle with an extra edge hanging off (see
  Figure~\ref{fig:Theta} in Example~\ref{ex:Theta}). By the
  recurrence~\eqref{eq:kappa} in Corollary~\ref{cor:kappa} along
  with~\eqref{eq:Circle_n},
  \[
  \kappa(\Gamma)=\kappa(\Gamma_e')+\kappa(\Gamma_e'')=4+2=6\,.
  \]
  Representatives from the six distinct $\kappa$-classes of $\Gamma$ are
  shown in Figure~\ref{fig:transversal}. This particular transversal
  was chosen so that $v$ is a source, and so the $vw$-intervals
  can be identified immediately, and they are (from left-to-right)
  \[
  \{v,v',w',w\}\,,\qquad \{v,w\}\,,\qquad \{v,v',w',z,w\}\,,\qquad
  \{v,z,w\}\,,\qquad \{v,z,w\}\,,\qquad \{v,w\}\,.
  \]
  Note that the first two orientations in Figure~\ref{fig:transversal}
  are $O^a_\Gamma$ and $O^b_\Gamma$, and the third is
  $\kappa$-equivalent to $O^c_\Gamma$. Letting $P=v,v',w',w,v$ and
  $Q=v,z,w,v$ be the paths as defined in Example~\ref{ex:nu}, the pair
  $(\nu_P,\nu_Q)$ is a complete invariant of
  $\Acyc(\Gamma)/\simkappa$. For each of the six orientations, the
  values of $\nu_P(O_\Gamma)$ (top) and $\nu_Q(O_\Gamma)$ (bottom) are
  shown in Figure~\ref{fig:transversal}.
\begin{figure}
\begin{tikzpicture}[scale=0.5]
  \tikzstyle{VertexStyle}=[shape=circle,minimum size=2.5pt,inner sep=1.6pt,draw]
  \Vertex[style={fill=white}, x=1.5, y=4, L=\tiny {}]{v}
  \Vertex[style={fill=white}, x=4.5, y=4, L=\tiny {}]{w}
  \Vertex[style={fill=white}, x=2.25, y=3, L=\tiny {}]{vv}
  \Vertex[style={fill=white}, x=3.75, y=3, L=\tiny {}]{ww}
  \Vertex[style={fill=white}, x=3, y=1, L=\tiny {}]{z}
  \Edge[style = {post}](v)(w)
  \Edge[style = {post}](v)(vv)
  \Edge[style = {post}](vv)(ww)
  \Edge[style = {post}](ww)(w)
  \Edge[style = {post, bend right}](v)(z)
  \Edge[style = {post, bend left}](w)(z)
  \draw[fill=white] (1.5,4) node[label=above:\small $v$]{};
  \draw[fill=white] (4.5,4) node[label=above:\small $w$]{};
  \draw[fill=white] (3,2.65) node[label=above:\small $2$]{};
  \draw[fill=white] (3,1.25) node[label=above:\small $-1$]{};
\end{tikzpicture}
\hspace{0.15in}
\begin{tikzpicture}[scale=0.5]
  \tikzstyle{VertexStyle}=[shape=circle,minimum size=2.5pt,inner sep=1.6pt,draw]
  \Vertex[style={fill=white}, x=1.5, y=4, L=\tiny {}]{v}
  \Vertex[style={fill=white}, x=4.5, y=4, L=\tiny {}]{w}
  \Vertex[style={fill=white}, x=2.25, y=3, L=\tiny {}]{vv}
  \Vertex[style={fill=white}, x=3.75, y=3, L=\tiny {}]{ww}
  \Vertex[style={fill=white}, x=3, y=1, L=\tiny {}]{z}
  \tikzstyle{LabelStyle}=[above]
  \Edge[style = {post}](v)(w)
  \Edge[style = {post}](v)(vv)
  \Edge[style = {pre}](vv)(ww)
  \Edge[style = {post}](ww)(w)
  \Edge[style = {post, bend right}](v)(z)
  \Edge[style = {post, bend left}](w)(z)
  \draw[fill=white] (1.5,4) node[label=above:\small $v$]{};
  \draw[fill=white] (4.5,4) node[label=above:\small $w$]{};
  \draw[fill=white] (3,2.65) node[label=above:\small $0$]{};
  \draw[fill=white] (3,1.25) node[label=above:\small $-1$]{};
\end{tikzpicture}
\hspace{0.15in}
\begin{tikzpicture}[scale=0.5]
  \tikzstyle{VertexStyle}=[shape=circle,minimum size=2.5pt,inner sep=1.6pt,draw]
  \Vertex[style={fill=white}, x=1.5, y=4, L=\tiny {}]{v}
  \Vertex[style={fill=white}, x=4.5, y=4, L=\tiny {}]{w}
  \Vertex[style={fill=white}, x=2.25, y=3, L=\tiny {}]{vv}
  \Vertex[style={fill=white}, x=3.75, y=3, L=\tiny {}]{ww}
  \Vertex[style={fill=white}, x=3, y=1, L=\tiny {}]{z}
  \tikzstyle{LabelStyle}=[above]
  \Edge[style = {post}](v)(w)
  \Edge[style = {post}](v)(vv)
  \Edge[style = {post}](vv)(ww)
  \Edge[style = {post}](ww)(w)
  \Edge[style = {post, bend right}](v)(z)
  \Edge[style = {pre, bend left}](w)(z)
  \draw[fill=white] (1.5,4) node[label=above:\small $v$]{};
  \draw[fill=white] (4.5,4) node[label=above:\small $w$]{};
  \draw[fill=white] (3,2.65) node[label=above:\small $2$]{};
  \draw[fill=white] (3,1.25) node[label=above:\small $1$]{};
\end{tikzpicture}
\hspace{0.15in}
\begin{tikzpicture}[scale=0.5]
  \tikzstyle{VertexStyle}=[shape=circle,minimum size=2.5pt,inner sep=1.6pt,draw]
  \Vertex[style={fill=white}, x=1.5, y=4, L=\tiny {}]{v}
  \Vertex[style={fill=white}, x=4.5, y=4, L=\tiny {}]{w}
  \Vertex[style={fill=white}, x=2.25, y=3, L=\tiny {}]{vv}
  \Vertex[style={fill=white}, x=3.75, y=3, L=\tiny {}]{ww}
  \Vertex[style={fill=white}, x=3, y=1, L=\tiny {}]{z}
  \tikzstyle{LabelStyle}=[above]
  \Edge[style = {post}](v)(w)
  \Edge[style = {post}](v)(vv)
  \Edge[style = {post}](vv)(ww)
  \Edge[style = {pre}](ww)(w)
  \Edge[style = {post, bend right}](v)(z)
  \Edge[style = {pre, bend left}](w)(z)
  \draw[fill=white] (1.5,4) node[label=above:\small $v$]{};
  \draw[fill=white] (4.5,4) node[label=above:\small $w$]{};
  \draw[fill=white] (3,2.65) node[label=above:\small $0$]{};
  \draw[fill=white] (3,1.25) node[label=above:\small $1$]{};
\end{tikzpicture}
\hspace{0.15in}
\begin{tikzpicture}[scale=0.5]
  \tikzstyle{VertexStyle}=[shape=circle,minimum size=2.5pt,inner sep=1.6pt,draw]
  \Vertex[style={fill=white}, x=1.5, y=4, L=\tiny {}]{v}
  \Vertex[style={fill=white}, x=4.5, y=4, L=\tiny {}]{w}
  \Vertex[style={fill=white}, x=2.25, y=3, L=\tiny {}]{vv}
  \Vertex[style={fill=white}, x=3.75, y=3, L=\tiny {}]{ww}
  \Vertex[style={fill=white}, x=3, y=1, L=\tiny {}]{z}
  \tikzstyle{LabelStyle}=[above]
  \Edge[style = {post}](v)(w)
  \Edge[style = {post}](v)(vv)
  \Edge[style = {pre}](vv)(ww)
  \Edge[style = {pre}](ww)(w)
  \Edge[style = {post, bend right}](v)(z)
  \Edge[style = {pre, bend left}](w)(z)
  \draw[fill=white] (1.5,4) node[label=above:\small $v$]{};
  \draw[fill=white] (4.5,4) node[label=above:\small $w$]{};
  \draw[fill=white] (3,2.65) node[label=above:\small $-2$]{};
  \draw[fill=white] (3,1.25) node[label=above:\small $1$]{};
\end{tikzpicture}
\hspace{0.15in}
\begin{tikzpicture}[scale=0.5]
  \tikzstyle{VertexStyle}=[shape=circle,minimum size=2.5pt,inner sep=1.6pt,draw]
  \Vertex[style={fill=white}, x=1.5, y=4, L=\tiny {}]{v}
  \Vertex[style={fill=white}, x=4.5, y=4, L=\tiny {}]{w}
  \Vertex[style={fill=white}, x=2.25, y=3, L=\tiny {}]{vv}
  \Vertex[style={fill=white}, x=3.75, y=3, L=\tiny {}]{ww}
  \Vertex[style={fill=white}, x=3, y=1, L=\tiny {}]{z}
  \tikzstyle{LabelStyle}=[above]
  \Edge[style = {post}](v)(w)
  \Edge[style = {post}](v)(vv)
  \Edge[style = {pre}](vv)(ww)
  \Edge[style = {pre}](ww)(w)
  \Edge[style = {post, bend right}](v)(z)
  \Edge[style = {post, bend left}](w)(z)
  \draw[fill=white] (1.5,4) node[label=above:\small $v$]{};
  \draw[fill=white] (4.5,4) node[label=above:\small $w$]{};
  \draw[fill=white] (3,2.65) node[label=above:\small $-2$]{};
  \draw[fill=white] (3,1.25) node[label=above:\small $-1$]{};
\end{tikzpicture}
\caption{A transversal for $\Acyc(\Gamma)/\simkappa$, with the values
  of the complete invariant $(\nu_P,\nu_Q)$.}
\label{fig:transversal}
\end{figure}
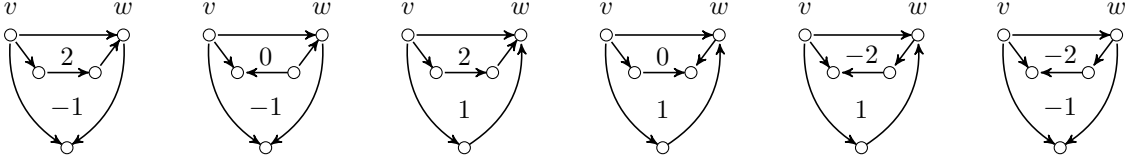
\end{ex}


\section{Conjugacy of Coxeter Elements.}
\label{sec:conjugacy}

Our analysis of $\Acyc(\Gamma)/\simkappa$ also gives a straightforward
solution to the conjugacy problem for Coxeter elements in simply-laced
Coxeter groups. Before stating the theorem and proof, we will briefly
review the connection between $\kappa$-equivalence and Coxeter theory,
as described in~\cite{Macauley:08b}. A \emph{Coxeter group} is a
generalized reflection group, generated by $n$ distinguished
involutions $s_1,\dots,s_n$ by the presentation
\[
W=\<s_1,\dots,s_n\mid (s_is_j)^{m_{ij}}\>\,,
\]
where $m_{ij}=1$ iff $i=j$, and $m_{ij}\geq 2$ otherwise. If $s_is_j$
has infinite order, then we say that $m_{ij}=\infty$. The pair $(W,S)$
of the group $W$ with the generating set $S$ is called a \emph{Coxeter
  system}, which is uniquely encoded by its Coxeter graph $\Gamma$,
with vertex set $S$ and edge set $\{s_i,s_j\}$ for which $m_{ij}\geq
3$, where each edge has weight $m_{ij}$. A Coxeter group is
\emph{simply-laced} if each $m_{ij}\leq 3$, thus every undirected
graph is the Coxeter graph of a simply-laced Coxeter group.

A \emph{Coxeter element} is the product of the generators in some
order, and there is a natural bijection between the set $\C(W,S)$ of
Coxeter elements $\prod_i s_{\pi(i)}$ of a Coxeter group (see,
\eg~\cite{Bjorner:05}) with generators $s_i$ for $1\le i\le n$ and
Coxeter graph $\Gamma$ (ignoring bond strengths), and the set of
$\alpha$-equivalence classes $[\pi]_\Gamma$. This is clear since the
commuting generators are precisely those that are not connected in
$\Gamma$. Thus, there is a natural bijection
\begin{equation}
  \label{eq:bij2}
  \C(W,S) \longrightarrow \Acyc(\Gamma) \,.
\end{equation}
Moreover, conjugating a Coxeter element $c=\prod s_{\pi(i)}$ by
$s_{\pi(1)}$ corresponds to a cyclic shift, \ie,
\begin{equation*}
  s_{\pi(1)} (s_{\pi(1)} s_{\pi(2)} \cdots s_{\pi(n)}) s_{\pi(1)}
  = s_{\pi(2)} \cdots s_{\pi(n)} s_{\pi(1)}\,,
\end{equation*}
since each generator $s_i$ is an involution. Therefore,
$\kappa$-equivalence naturally carries over to an equivalence relation
on $\C(W,S)$. The $\nu$-function carries over as well -- define
$\nu_P(c)$ to be $\nu_P(O_\Gamma)$, where $O_\Gamma$ is the acyclic
orientation of $\Gamma$ corresponding to $c$. It is now elementary to
see that $c,c'\in\C(W,S)$ are conjugate if $c\tsimkappa c'$. However,
the converse of this statement was not proven until 2009.
\begin{thm}[\cite{Eriksson:09}]
  \label{thm:conjugacy}
  Let $(W,S)$ be a Coxeter system. Then two Coxeter elements $c,c'\in
  \C(W,S)$ are conjugate if and only if $c\tsimkappa c'$.
\end{thm}
It follows immediately that the number of distinct conjugacy classes
containing Coxeter elements is exactly $\kappa(\Gamma)$. Until the
Erikssons' proof, the result was known only for the special case of
$\C(W,S)$ when $\Gamma$ was simply-laced and unicyclic, established by
Shi in 2001~\cite{Shi:01}. It is elementary to weaken the simply-laced
condition to the bond strengths being multiples of three or infinite,
which Shi mentions in~\cite{Shi:01}. The bijection in
Theorem~\ref{thm:bijection} applied to Shi's result for unicyclic
Coxeter graphs yields a simple and elegant proof of the result for all
simply-laced systems, which we present below.

\begin{proof}
  Suppose for sake of contradiction that the Coxeter elements
  $c=c_1c_2\cdots c_n$ and $c'=c'_1c'_2\cdots c'_n$ are conjugate with
  $c\nkeq c'$, and that $wcw^{-1} = c'$ for some $w=w_1\dots w_k\in W$
  with each $w_i\in S$. By Theorem~\ref{thm:nu}, there is some simple
  chordless cycle $P=v_0,v_1,\dots,v_{m-1},v_m$ (i.e., $v_0=v_m$) in
  $\Gamma$ such that $\nu_P(c)\neq\nu_P(c')$.  Let $S_P=\{s_i\mid i\in
  P\}$, and let $C_m$ be the (circular) Coxeter graph induced by the
  vertices in $P$. The Coxeter group generated by $S_P$ is the affine
  Weyl group $\widetilde{A}_{m-1}$, and there is a natural
  homomorphism
  $W\stackrel{\varphi}{\longrightarrow}\widetilde{A}_{m-1}$ defined on
  the generators by
  \[
  s_i\stackrel{\varphi}{\longmapsto}
  \left\{\begin{array}{cl} s_i & \quad i\in P \\ 1 & \quad i\not\in
  P\,.\end{array}\right.
  \]
  Since $c$ and $c'$ are conjugate in $W$, $\varphi(c)$ and
  $\varphi(c')$ are conjugate in $\widetilde{A}_{m-1}$. By choice of
  $P$, $\nu_P(c)\neq\nu_P(c')$, and thus
  $\varphi(c)\nkeq\varphi(c')$. However, since the statement holds for
  unicyclic graphs, we must have $\varphi(c)\tsimkappa\varphi(c')$,
  which is the desired contradiction.
\end{proof}
Theorems~\ref{thm:nu} and~\ref{thm:conjugacy} give us an easy way to
verify whether two Coxeter elements $c$ and $c'$ in any Coxeter group
$W$ are conjugate. Pick a cycle basis of the Coxeter graph $\Gamma$,
and for each cycle $C$, compute $\nu_C(c)$ and $\nu_C(c')$. By
Theorem~\ref{thm:nu}, $c\tsimkappa c'$ iff $\nu_C(c)=\nu_C(c')$ for
each $C$. By Theorem~\ref{thm:conjugacy}, this is equivalent to $c$
and $c'$ being conjugate in $W$. Therefore, conjugacy of Coxeter
elements can be verified in $O(n^2)$ steps, where
$n=|\vset[\Gamma]|$. One application of this is seeing how the
conjugacy classes split as an edge $\{s_i,s_j\}$ is added to $\Gamma$
(or equivalently, as the relation $(s_is_j)^{m_{ij}}$ is added to the
group presentation).


\section{Discrete Dynamical Systems, Node-firing Games, and Quiver
  Representations.}
\label{sec:summary}

We conclude with a brief discussion of how the equivalence relation
studied in this paper arises in various areas of mathematics. The
original motivation came from both authors' interest in
\emph{sequential dynamical systems} (\sdss). The equivalence relation
$\tsimalpha$ arises naturally in the study of functional equivalence
of these systems. This can be seen as follows. Given a graph $\Gamma$
with vertex set $\{1,2,\ldots,n\}$ as above, a state $x_v\in K$ is
assigned to each vertex $v$ of $\Gamma$ for some finite set $K$.  The
\emph{system state} is the tuple consisting of all the vertex states,
and is denoted by $x = (x_1,\ldots,x_n)\in K^n$. The sequence of states
associated to the $1$-neighborhood $B_1(v;\Gamma)$ of $v$ in $\Gamma$
(in some fixed order) is denoted by $x[v]$. A sequence of vertex
functions $(f_i)_i$ with $f_i \colon K^{d(i)+1} \longrightarrow K$
induces $\Gamma$-local functions $F_i \colon K^n \longrightarrow K^n$
of the form
\begin{equation*}
  F_i(x_1,\ldots,x_n) =
  (x_1,\ldots,x_{i-1}, f_i(x[i]), x_{i+1}, \ldots, x_n)\,.
\end{equation*}
The sequential dynamical system map with update order $\pi = (\pi_i)_i
\in S_\Gamma$ is the function composition
\begin{equation}
  F_\pi = F_{\pi_n} \circ F_{\pi_{n-1}} \circ \cdots \circ
  F_{\pi_2} \circ F_{\pi_1} \,.
\end{equation}
By construction, if $\pi\tsimalpha\pi'$ holds, then $F_\pi$ and
$F_{\pi'}$ are identical as functions, independent of the choice of
state space $K$ or vertex functions. Thus, $\acyc(\Gamma)$ is a
general upper bound for the number of functionally non-equivalent \sds
maps that can be generated over the graph $\Gamma$ for a fixed
sequence of $\Gamma$-local functions. Moreover, for any graph
$\Gamma$, there exist $\Gamma$-local functions for which this bound is
sharp~\cite{Mortveit:01a}.
A weaker form of equivalence is \textit{cycle equivalence}, which
means that the dynamical system maps are conjugate (using the discrete
topology) when restricted to their sets of periodic points. In the
language of graph theory, this means their periodic orbits are
isomorphic as directed graphs. For an update order
$\pi=\pi_1\cdots\pi_n$, define
$\operatorname{shift}(\pi)=\pi_2\cdots\pi_n\pi_1$.  The following
result shows how $\kappa$-equivalent update orders yield dynamical
system maps that are cycle equivalent.
\begin{thm}
  \label{thm:shift}
  For any finite set $K$ of vertex states, and for any $\pi\in
  S_\Gamma$, the \sds maps $F_\pi$ and
  $F_{\operatorname{shift}(\pi)}$ are cycle equivalent.
\end{thm}
We refer to~\cite{Macauley:09a} for the proof of this result, as well
as additional background on equivalences of sequential dynamical
systems, and applications of $\kappa$-equivalence to the structural
properties of their phase spaces. It is interesting to note that for
the class of \emph{update sequence independent}
(see~\cite{Macauley:11a}) sequential dynamical systems with binary
states, there is an additional equivalence on acyclic orientations that
governs cycle equivalence: reversal of all edge orientations.

\medskip

The \textit{chip-firing} game was introduced by Bj\"{o}rner,
Lov\'{a}sz, and Shor~\cite{Bjorner:91}. It is played over an
undirected graph $\Gamma$, and each vertex is given some number of
(but possibly zero) chips. If vertex $i$ has degree $d_i$, and at
least $d_i$ chips, then a legal move (or a ``click'') of vertex $i$ is
a transfer of one chip to each neighboring vertex. This may be viewed
as a generalization of a source-to-sink move for acyclic orientations
where the out-degree of a vertex plays the role of the chip count. The
chip-firing game is closely related to the \textit{numbers
  game}~\cite{Bjorner:05}. In the numbers game over a graph $\Gamma$,
the legal sequences of moves are in 1--1 correspondence with the
reduced words of the Coxeter group with Coxeter graph $\Gamma$. For an
excellent summary and comparison of these games,
see~\cite{KEriksson:94}.

\medskip

A quiver is a finite directed graph (loops and multiple edges are
allowed), and appears primarily in the study of representation
theory. A quiver $Q$ with a field $K$ gives rise to a path algebra
$KQ$, and there is a natural correspondence between $KQ$-modules and
$K$-representations of $Q$. In fact, there is an equivalence between
the categories of quiver representations, and modules over path
algebras. A path algebra is finite-dimensional if and only if the
quiver is acyclic, and the modules over finite-dimensional path
algebras form a reflective subcategory. A \textit{reflection functor}
maps representations of a quiver $Q$ to representations of a quiver
$Q'$, where $Q'$ differs from $Q$ by a source-to-sink
operation~\cite{Marsh:03}. We note that while the composition of $n$
source-to-sink operations (one for each vertex) maps a quiver back to
itself, the corresponding composition of reflection functors is not
the identity, but rather a \textit{Coxeter functor}. In fact, the same
result in~\cite{Speyer:09} about powers of Coxeter elements being
reduced was proven previously using techniques from the representation
theory of quivers~\cite{Kleiner:07}.

\medskip

We hope this paper will motivate further explorations of the
connections between these topics. We are particularly curious about
any implications to the representation theory of quivers. This is a
field which the both authors of this paper are quite unfamiliar with,
yet it motivated Kleiner and Pelley to study admissible sequences and
apply these tools from quiver representations to Coxeter
groups. Without this work, the aforementioned papers of Speyer and the
Erikssons would likely not have materialized.

\bigskip


\noindent\textbf{Acknowledgments.} Both authors are grateful to the
NDSSL group at Virginia Tech for the support of this research. Special
thanks to Ed Green and Ken Millett for helpful discussions and
feedback, and to William Y.~C. Chen for valuable advice regarding the
preparation and structuring of this paper. 



\begin{thebibliography}{10}

\bibitem{Auslander:97}
M.~Auslander, I.~Reiten, and S.~O. Smal\o.
\newblock {\em Representation Theory of Artin Algebras}.
\newblock Cambridge University Press, 1997.

\bibitem{Bjorner:05}
A.~Bj{\"o}rner and F.~Brenti.
\newblock {\em Combinatorics of {C}oxeter Groups}.
\newblock Springer-Verlag, New York, 2005.

\bibitem{Bjorner:91}
A.~Bj\"{o}rner, L.~Lov\'{a}sz, and P.~Shor.
\newblock Chip-firing games on graphs.
\newblock {\em European J. Combin.}, 12:283--291, 1991.

\bibitem{Cartier:69}
P.~Cartier and D.~Foata.
\newblock {\em Problemes combinatoires de commutation et
  re\'arrangements},
  volume~85 of {\em Lect. Notes Math.}.
\newblock Springer Verlag, 1969.

\bibitem{Coleman:89}
A.~J. Coleman.
\newblock Killing and the {C}oxeter transformation of {K}ac-{M}oody
algebras.
\newblock {\em Invent. Math.}, 95:447--477, 1989.

\bibitem{Eriksson:09}
H.~Eriksson and K.~Eriksson.
\newblock Conjugacy of {C}oxeter elements.
\newblock {\em Electron. J. Combin.}, 16(2):\#R4, 2009.

\bibitem{KEriksson:94}
K.~Eriksson.
\newblock Node firing games on graphs.
\newblock {\em Contemp. Math.}, 178:117--127, 1994.

\bibitem{Kleiner:07}
M.~Kleiner and A.~Pelley.
\newblock Admissible sequences, preprojective representations of
quivers, and
  reduced words in the {W}eyl group of a {K}ac-{M}oody algebra.
\newblock {\em Internat. Math. Res. Notices}, 2007, May
2007.

\bibitem{Macauley:11a}
M.~Macauley, J.~McCammond, and H.~S. Mortveit.
\newblock Dynamics groups of asynchronous cellular automata.
\newblock {\em J. Algebraic Combin.}, 33:31--55, 2011.

\bibitem{Macauley:08b}
M.~Macauley and H.~S. Mortveit.
\newblock On enumeration of conjugacy classes of {C}oxeter elements.
\newblock {\em Proc. Amer. Math. Soc.}, 136:4157--4165, 2008.

\bibitem{Macauley:09a}
M.~Macauley and H.~S. Mortveit.
\newblock Cycle equivalence of graph dynamical systems.
\newblock {\em Nonlinearity}, 22:421--436, 2009.

\bibitem{Marsh:03}
R.~Marsh, M.~Reineke, and A.~Zelevinsky.
\newblock Generalized associahedra via quiver representations.
\newblock {\em Trans. Amer. Math. Soc.},
  355:4171--4186, 2003.

\bibitem{Miller:07}
E.~Miller, V.~Reiner, and B.~Sturmfels, editors.
\newblock {\em Geometric Combinatorics}, volume~13 of
          {\em IAS/Park City Mathematics Series}.
\newblock AMS and IAS/Park City Mathematics Institute, 2007.

\bibitem{Mortveit:01a}
H.~S. Mortveit and C.~M. Reidys.
\newblock Discrete, sequential dynamical systems.
\newblock {\em Discrete Math.}, 226:281--295, 2001.

\bibitem{Novik:02}
I.~Novik, A.~Postnikov, and B.~Sturmfels.
\newblock Syzygies of oriented matroids.
\newblock {\em Duke Math. J.}, 111:287--317, 2002.

\bibitem{Pretzel:86}
O.~Pretzel.
\newblock On reorienting graphs by pushing down maximal vertices.
\newblock {\em Order}, 3(2):135--153, 1986.

\bibitem{Shi:97a}
J.-Y. Shi.
\newblock The enumeration of {C}oxeter elements.
\newblock {\em J. Algebraic Combin.}, 6:161--171, 1997.

\bibitem{Shi:01}
J.-Y. Shi.
\newblock Conjugacy relation on {C}oxeter elements.
\newblock {\em Adv. Math.}, 161:1--19, 2001.

\bibitem{Speyer:09}
D.~E. Speyer.
\newblock Powers of {C}oxeter elements in infinite groups are reduced.
\newblock {\em Proc. Amer. Math. Soc.}, 137:1295--1302, 2009.

\bibitem{Stembridge:96}
J.~R. Stembridge.
\newblock On the fully commutative elements of {C}oxeter groups.
\newblock {\em J. Algebraic Combin.}, 5:353--385, 1996.

\end{thebibliography}
\end{document}